\documentclass[12pt,a4paper,draft]{article}
\usepackage[english,russian]{babel}
\usepackage{amssymb, amsmath, latexsym, amsfonts, amsthm}
\begin{document}
\renewcommand{\refname}{References}
\newtheorem{theorem}{Theorem}
\newtheorem{lemma}{Lemma}
\newtheorem{corollary}{Corollary}
\newtheorem{proposition}{Proposition}
\begin{center}
On a space of entire functions and its Fourier transform
\end{center}
\begin{center}
I.Kh. Musin\footnote {E-mail: musin\_ildar@mail.ru}
\end{center}
\begin{center}
Institute of Mathematics with Computer Centre of Ufa Scientific Centre of Russian Academy of Sciences, 

Chernyshevsky str., 112, Ufa, 450077, Russia
\end{center}

\vspace {0.3cm}

\renewcommand{\abstractname}{}
\begin{abstract}
{\sc Abstract}. A space of entire functions of several complex variables rapidly decreasing on ${\mathbb R}^n$ and such that their growth along $i{\mathbb R}^n$ is majorized with the help of a family of weight functions is considered in this paper. For such space an equivalent description in terms of estimates on all of its partial derivatives as functions on ${\mathbb R}^n$ and a Paley-Wiener type theorem are obtained. 

\vspace {0.3cm}
MSC: 32A15, 42B10, 46E10, 46F05, 42A38

%46F05   	Topological linear spaces of test functions, distributions and ultradistributions

%46E10   	Topological linear spaces of continuous, differentiable or analytic functions

%32A15   	Entire functions

%42B10   	Fourier and Fourier-Stieltjes transforms and other transforms of Fourier type (for n >1)

%42A38   	Fourier and Fourier-Stieltjes transforms and other transforms of Fourier type (for n=1)

\vspace {0.3cm}
Keywords: Gelfand-Shilov spaces, Fourier transform, entire functions, convex functions 

\end{abstract}

\vspace {1cm}
 
\section{Introduction}

{\bf 1.1. On the problem}. In the 1950's the study of $W$-type spaces started  with the works of B.L. Gurevich \cite {Gur1}, \cite {Gur2} and I.M. Gelfand and G.E. Shilov \cite {GS1}, \cite {GS2}. They described them by means of the Fourier transformat and then applied this description to study the uniqueness of the Cauchy problem of partial differential equations and their systems. 
These spaces generalize Gelfand-Shilov spaces of $S$-type  \cite {GS1}. So they are often called  Gelfand-Shilov spaces of $W$-type. 

Let us recall the definition of the Gelfand-Shilov spaces of $W$-type. 
Let $M$ and $\Omega$ be differentiable functions on $[0, \infty)$ such that 
$M(0)= \Omega(0) = M'(0) = \Omega'(0) = 0$ and, moreover, so that their derivatives are continuous, increasing and unbounded at infinity. 
Considering ${\mathbb R}^n$ with its usual norm $\Vert \cdot \Vert$, 
$W_M$ is the space of all infinitely differentiable functions
$f$ on ${\mathbb R}^n$ satisfying the following upper estimate on every partial derivative
$$
\vert (D^{\alpha} f) (x) \vert \le C_{\alpha} e^{-M(a\Vert x \Vert)}, \ x \in {\mathbb R}^n,
$$
for some positive constant $a$.
Also, $W^{\Omega}$ is the space of entire functions $f$ on ${\mathbb C}^n$ satisfying the estimate 
$$
\vert \zeta^{\alpha} f(\zeta) \vert \le C_{\alpha} 
e^{\Omega(b \Vert \eta \Vert)}, \ \zeta = \xi +i \eta \in {\mathbb R}^n + i{\mathbb R}^n, \alpha  \in {\mathbb Z}_+^n,
$$ 
for some $b > 0$.
And finally, $W_M^{\Omega}$ is the space of entire functions $f$ on ${\mathbb C}^n$ satisfying
the estimate 
$$
\vert f(\xi +i \eta) \vert \le C e^{-M(a\Vert \xi \Vert) + \Omega(b \Vert \eta \Vert)}, \ \xi, \eta \in {\mathbb R}^n,
$$ 
for some positive constants $a, b$ and $C$.

$W$-type spaces and some their generalizations have been studied by many mathematicians. For example, new characterizations of $W$-type spaces and their generalizations were given  by J. Chung, S.Y. Chung and D. Kim \cite {C-C-K 1}, \cite {C-C-K 2},
R.S. Pathak and S.K. Upadhyay \cite {P-U}, S.K. Upadhyay \cite {U} 
(in terms of Fourier transform), 
N.G. De Bruijn \cite {Br}, A.J.E.M. Janssen and S.J.L. van Eijndhoven \cite {J-E}, Jonggyu Cho \cite {Cho}
(by using the growth of their Wigner distributions). R.S. Pathak \cite {P} and S.J.L. van Eijndhoven and M.J. Kerkhof \cite {E-K} introduced new spaces of W-type and investigated the behaviour of the Hankel transform over
them (see also  \cite {B-R}). New $W$-type spaces introduced by V.Ya. Fainberg, M.A. Soloviev \cite {FS} turned out to be useful for nonlocal theory of highly singular quantum fields. 

In this paper we explore spaces of entire functions which are 
natural generalizations of $W^{\Omega}$-type spaces.
Namely, we work with the space $E(\varPhi)$ of entire functions that we introduce next.
Let $n \in {\mathbb N}$, $H({\mathbb C}^n)$ be the space of entire functions on ${\mathbb C}^n$,  
$\Vert u \Vert$ be the Euclidean norm of $u \in {\mathbb R}^n ({\mathbb C}^n)$. 
Denote by ${\mathcal A}({\mathbb R}^n)$ the set of all real-valued functions $g \in C({\mathbb R}^n)$ satisfying the following  conditions:

1) $g(x) = g(\vert x_1 \vert, \ldots , \vert x_n \vert), \ x = (x_1, \ldots , x_n)\in {\mathbb R}^n$;

2) the restriction of $g$ to $[0, \infty)^n$ is nondecreasing in each variable;
%$g(x_1, \ldots , x_n)$ is nondecreasing in each variable in $[0, \infty)^n$;

3) $\displaystyle \lim_{x \to \infty} \frac {g(x)}{\Vert x \Vert}= + \infty$.

\noindent Let $\varPhi =\{\varphi_{\nu}\}_{{\nu}=1}^{\infty}$ be a subset of ${\mathcal A}({\mathbb R}^n)$ consisting of 
functions $\varphi_{\nu}$ satisfying the condition:

$i_0$) for each $\nu \in {\mathbb N}$ and each $A > 0$ there exists a constant $C_{\nu, A} > 0$ such that 
$$
\varphi_{\nu}(x) + A \ln (1 + \Vert x \Vert) \le \varphi_{\nu+1}(x) + C_{{\nu}, A}, \ x \in {\mathbb R}^n.
$$
\noindent For each $\nu \in {\mathbb N}$ and $m \in {\mathbb Z}_+$ consider the normed space
$$
E_m(\varphi_{\nu}) = \{f \in H({\mathbb C}^n): p_{\nu, m}(f) = \sup_{z \in {\mathbb C}^n} 
\frac 
{\vert f(z)\vert (1 + \Vert z \Vert)^m}
{e^{\varphi_{\nu} (Im \, z)}} < \infty \}.
$$
Let $E(\varphi_{\nu})= \bigcap \limits_{m=0}^{\infty} E_m(\varphi_{\nu})$.
%, $E(\varPhi)= \bigcup \limits_{\nu=1}^{\infty} E(\varphi_{\nu})$. 
%Since $p_{\nu, k}(f) \le p_{\nu, k+1}(f)$ 
%for $f \in E_{k+1}(\varphi_{\nu})$ 
%then 
Obviously, $E_{m+1}(\varphi_{\nu})$ is continuously embedded in $E_m(\varphi_{\nu})$. 
Endow $E(\varphi_{\nu})$ with a projective limit topology of spaces $E_m(\varphi_{\nu})$. 
Note that if $f \in E(\varphi_{\nu})$ then 
$p_{\nu+1, m}(f) \le e^{C_{\nu, 1}} p_{\nu, m}(f)$ for each $m \in {\mathbb Z}_+$. 
Hence, $E(\varphi_{\nu})$ is continuously embedded in $E(\varphi_{\nu + 1})$ for each $\nu \in {\mathbb N}$. Let $E(\varPhi)= \bigcup \limits_{\nu=1}^{\infty} E(\varphi_{\nu})$. With the usual operations of addition and multiplication by complex numbers 
$E(\varPhi)$
%and $E(\varPhi)$ are 
is a linear space. 
Supply $E(\varPhi)$ with a topology of the inductive limit of spaces $E(\varphi_{\nu})$.

In the paper we describe the space $E(\varPhi)$ in terms of estimates 
on partial derivatives of functions on ${\mathbb R}^n$ and study the Fourier transform of functions of $E(\varPhi)$ under additional conditions on $\varPhi$. Results of the paper could be useful in harmonic analysis, theory of entire functions, convex analysis and in the study of partial differential and pseudo-differential operators.

{\bf 1.2. Some notations}.  
For $u=(u_1, \ldots , u_n), v=(v_1, \ldots , v_n) \in {\mathbb R}^n \ ({\mathbb C}^n)$ let denote
$\langle u, v \rangle  = u_1 v_1 + \cdots + u_n v_n$. 

For $\alpha = (\alpha_1, \ldots , \alpha_n) \in {\mathbb Z}_+^n$, 
$x =(x_1, \ldots , x_n) \in {\mathbb R}^n$, 
$z =(z_1, \ldots , z_n) \in {\mathbb C}^n$ we follow the standard notations $\vert \alpha \vert = \alpha_1 + \ldots  + \alpha_n$, 
$\alpha! = \alpha_1! \cdots \alpha_n!$, 
$x^{\alpha} = x_1^{\alpha_1} \cdots x_n^{\alpha_n}$, $z^{\alpha} = z_1^{\alpha_1} \cdots z_n^{\alpha_n}$, 
$D^{\alpha}=
\frac {{\partial}^{\vert \alpha \vert}}{\partial x_1^{\alpha_1} \cdots \partial x_n^{\alpha_n}}$ .

%If the first $k$ ($1 \le k \le n-1$) coordinates of $x = (x_1, \ldots , x_n) \in {\mathbb R}^n$ are non-zero and other %coordinates are equal to zero then put $x' = (x_1, \ldots , x_k)$, 
%${\bf 0} = (0, \ldots, 0) \in {\mathbb R}^{n-k}$, $x = (x', {\bf 0})$.

By $s_n$ denote the surface area of the unit sphere in ${\mathbb R}^n$.
 
%For $r > 0$ let $\Pi_r: = \{x = (x_1, \ldots , x_n) \in {\mathbb R}^n: x_1 > r, \ldots , x_n > r \}$.

%If $X \subset {\mathbb R}^n$ then $int \, X$ is its interior.
%$\overline X$ is its closure
%in ${\mathbb R}^n$ and  

%Sometimes it will be convenient to denote $[0, \infty)^n$ by $[0, \infty)^n$.

If $[0, \infty)^n \subseteq X \subset {\mathbb R}^n$ then for a function $u$ on $X$ denote by $u[e]$ the function defined by the rule:
$u[e](x) = u(e^{x_1}, \ldots, e^{x_n}), \ x = (x_1, \ldots , x_n) \in  {\mathbb R}^n$.
%, and by $R(u)$ denote the restriction of $u$ to $[0, \infty)^n$.

%For brevity denote $\varphi_m [e]$ by $\psi_m$.

For an unbounded subset $X$ of ${\mathbb R}^n$ we denote by ${\mathcal B}(X)$ the set of all real-valued continuous functions $g$ on $X$ such that
$\displaystyle \lim_{x \to \infty, \atop x \in X} \frac {g(x)}{\Vert x \Vert}= + \infty $.

%Let
%${\mathcal B}_+ =\{g \in {\mathcal B}({\mathbb R}^n): 
%g(x) = g(\vert x_1 \vert, \ldots , \vert x_n \vert), \ x = (x_1, \ldots , x_n) \in {\mathbb R}^n \}$.

%Let $V =\{g \in {\mathcal B}({\mathbb R}^n): \text {$g$ is convex on ${\mathbb R}^n$} \}$.

Let us recall that the Young-Fenchel conjugate of a function $g:{\mathbb R}^n \to [-\infty, + \infty]$ is the function 
$g^*:{\mathbb R}^n \to [-\infty, + \infty]$ defined by
$$
g^*(x) = \displaystyle \sup \limits_{y \in {\mathbb R}^n}(\langle x, y \rangle - g(y)), \ x \in {\mathbb R}^n. 
$$

%let $dom \, g =\{x \in {\mathbb R}^n: g(x) < \infty \}$ and 

%%%%% Recall that if $g \in V$ then by the inversion formula for Young transform \cite {R} $(g^*)^*=g$.

{\bf 1.3. Main results and organization of the paper}. Given a family $\varPhi =\{\varphi_{\nu}\}_{\nu=1}^{\infty}$
as before and denoting $\varphi_{\nu} [e]$ by $\psi_{\nu}$, consider the family $\Psi^*=\{\psi_{\nu}^*\}_{\nu=1}^{\infty}$. 
For each $\nu \in {\mathbb N}$ and $m \in {\mathbb Z}_+$ consider the normed space 
$$
{\mathcal E}_m(\psi_{\nu}^*) =\{f \in  C^{\infty}({\mathbb R}^n): 
\rho_{m, \nu}(f)= \sup_{x \in {\mathbb R}^n, \alpha \in {\mathbb Z}_+^n} 
\frac {(1+ \Vert x \Vert)^m \vert (D^{\alpha}f)(x) \vert}{\alpha! e^{-\psi_{\nu}^*(\alpha)}} < \infty \}.
$$
Let ${\mathcal E}(\psi_{\nu}^*) = \bigcap \limits_{m=0}^{\infty}{\mathcal E}_m(\psi_{\nu}^*)$.
Endow ${\mathcal E}(\psi_{\nu}^*)$ with the topology defined by the family of norms $\rho_{m, \nu}$ ($m \in {\mathbb Z}_+$). 
Let ${\mathcal E}(\Psi^*) = \bigcup \limits_{\nu=1}^{\infty}{\mathcal E}(\psi_{\nu}^*)$. 
Supply ${\mathcal E}(\Psi^*)$ with an inductive limit topology of spaces ${\mathcal E}(\psi_{\nu}^*)$. 

The first two theorems are aimed to characterize functions of the space $E(\varPhi)$ 
in terms of estimates of their partial derivatives on ${\mathbb R}^n$.   

\begin{theorem} 
For each $f \in  E(\varPhi)$ its restriction to ${\mathbb R}^n$
%, $f_{|{\mathbb R}^n}$, 
belongs to ${\mathcal E}(\Psi^*)$.
\end{theorem}

\begin{theorem} 
Let the family $\varPhi=\{\varphi_{\nu}\}_{\nu=1}^{\infty}$ satisfies the additional conditions:

$i_1)$ for each $\nu \in {\mathbb N}$ there exist constants $\sigma_{\nu} > 1$ and $\gamma_{\nu} > 0$ such that
$$
\varphi_{\nu}(\sigma_{\nu} x) \le \varphi_{{\nu}+1}(x) + \gamma_{\nu}, \ x \in {\mathbb R}^n;
$$

$i_2)$ for each $\nu \in {\mathbb N}$ there exists a constant $K_{\nu} > 0$ such that 
$$
\varphi_{\nu}(x + \xi) \le \varphi_{{\nu}+1}(x) + K_{\nu}, \ x \in [0, \infty)^n, \xi \in [0, 1]^n.
$$

Then each function $f \in {\mathcal E}(\Psi^*)$ admits an (unique) extension to entire function belonging to $E(\varPhi)$.
\end{theorem}

The proofs of these theorems are given in section 3 by standard techniques. 
%Further, for $\nu \in {\mathbb N}$ and $k \in {\mathbb Z}_+$ define the space 
%$$
%{\mathcal H}_k(\varphi_{\nu})= \{f \in H({\mathbb C}^n): 
%$$
%$$
%{\cal N}_{\varphi_{\nu, k}}(f) = 
%\sup_{z \in {\mathbb C}^n} 
%\frac 
%{\vert f(z)\vert (1 + \Vert z \Vert)^k}
%{e^{(\psi_{\nu}^*)^*(\ln (1 + \vert Im z_1 \vert), \ldots , \ln (1 + \vert Im z_n \vert ))}} < \infty \}. 
%$$
%Let for each $\nu \in {\mathbb N}$ \ ${\mathcal H}(\varphi_{\nu})$ be a projective limit of spaces 
%${\mathcal H}_k(\varphi_{\nu})$. 
%Define the space ${\mathcal H}(\varPhi)$ as an inductive limit of spaces ${\mathcal H}(\varphi_{\nu})$. 
These proofs allow us to obtain additional information on the structure of the space $E(\varPhi)$ (see Proposition 2).

Section 4 is devoted to description of the space $E(\varPhi)$ in terms of Fourier transform under additional conditions on $\varPhi$. For each $\nu \in {\mathbb N}$ and $m \in {\mathbb Z}_+$ define the normed space
$$
G_m({\psi_{\nu}^*}) = \{f \in C^m({\mathbb R}^n): 
\Vert f \Vert_{m, \psi_{\nu}^*}
 = \sup_{x \in {\mathbb R}^n, \vert \alpha \vert \le m, \beta \in {\mathbb Z}_+^n}  
\frac 
{\vert x^{\beta}(D^{\alpha}f)(x) \vert}
{\beta! e^{-\psi_{\nu}^*(\beta)}} < \infty \}.
$$
Let $G({\psi_{\nu}^*})= \bigcap \limits_{m = 0}^{\infty} G_m({\psi_{\nu}^*})$.
%With usual operations of addition and multiplication by complex numbers $G({\psi_{\nu}^*})$ and $G({\Psi^*})$ are %linear spaces. 
Endow $G({\psi_{\nu}^*})$ with the topology defined by the family of norms 
$\Vert \cdot \Vert_{m, \psi_{\nu}^*}$ ($m \in {\mathbb Z}_+$).
Considering $G({\Psi^*})= \bigcup \limits_{\nu = 1}^{\infty} G({\psi_{\nu}^*})$ 
supply it with the topology of the inductive limit of spaces $G({\psi_{\nu}^*})$.

For $f \in E(\varPhi)$ define its Fourier transform $\hat f$ by the formula
$$
\hat f(x) = \int_{{\mathbb R}^n} f(\xi) e^{-i \langle x, \xi \rangle} \ d \xi , \ x \in {\mathbb R}^n.
$$

\begin{theorem} 
Let the family $\varPhi=\{\varphi_{\nu}\}_{{\nu}=1}^{\infty}$ satisfies  the condition $i_2)$ of Theorem 2 and the following two conditions:

$i_3)$ for each $\nu \in {\mathbb N}$ there is a constant $a_{\nu} > 0$ such that 
$$
\varphi_{\nu}(2x) \le \varphi_{\nu+1}(x) + a_{\nu}, \ x \in {\mathbb R}^n;
$$

$i_4)$ for each ${\nu} \in {\mathbb N}$ there is a constant $l_{\nu} > 0$ such that 
$$
2\varphi_{\nu}(x) \le \varphi_{\nu + 1}(x) + l_{\nu}, \ x \in {\mathbb R}^n.
$$
 
Then Fourier transform ${\mathcal F}: f \in E(\varPhi) \to  \hat f$ 
establishes an isomorphism between the spaces $E(\varPhi)$ and $G({\Psi^*})$.
\end{theorem}

Moreover, let ${\varPhi^*}=\{{\varphi_{\nu}^*}\}_{\nu=1}^{\infty}$. 
For each $\nu \in {\mathbb N}$ and $m \in {\mathbb Z}_+$ define the space
$$
GS_m(\varphi_{\nu}^*) = \{f \in C^m({\mathbb R}^n): q_{m, \nu}(f) = 
\sup_{x \in {\mathbb R}^n, \atop \vert \alpha \vert \le m} 
\frac {\vert (D^{\alpha}f)(x) \vert} {e^{-\varphi_{\nu}^*(x)}} < \infty \}.
$$
For each $\nu \in {\mathbb N}$ let
$GS({\varphi_{\nu}^*}) = \bigcap \limits_{m \in {\mathbb Z_+}}GS_m({\varphi_{\nu}^*})$. 
%With usual operations of addition and multiplication by complex numbers $GS({\varphi_{\nu}^*})$ and $GS({\varPhi^*})$ %are linear spaces.
%Note that for each $\nu \in {\mathbb N}$ and $m \in {\mathbb Z}_+$
%$$
%q_{m, \nu}(f) \le q_{m + 1, \nu}(f), \ f \in GS_{m+1}(\varphi_{\nu}^*);  
%$$
%$$
%q_{m, \nu + 1}(f) \le e^{C(\nu, 1)} q_{m, \nu}(f) , \ f \in GS_m(\varphi_{\nu}^*)).  
%$$
%Hence, $GS_{m+1}(\varphi_{\nu}^*)$ is continuously embedded in $GS_m(\varphi_{\nu}^*)$. 
Endow $GS({\varphi_{\nu}^*})$ with the topology defined by the family of norms $q_{m, \nu}$ ($m \in {\mathbb Z}_+$). 
%From the last inequality it follows that $GS({\varphi_{\nu}^*})$ 
%is continuously embedded in $GS({\varphi_{\nu + 1}^*})$. 
Let 
$GS({\varPhi^*})= \bigcup \limits_{\nu \in {\mathbb N}}GS({\varphi_{\nu}^*})$. Supply $GS({\varPhi^*})$ with an inductive limit topology of spaces $GS({\varphi_{\nu}^*})$. 

The main result of the section 5 is the following theorem.

\begin{theorem}
Let functions of the family $\varPhi$ be convex and satisfy the condition $i_2)$ of Theorem 2 and the condition 
$i_3)$ of Theorem 3. 
%for each $m \in {\mathbb N}$ there is a constant $a_m > 0$ such that 
%$$
%\varphi_m(2x) \le \varphi_{m+1}(x) + a_m, \ x \in {\mathbb R}^n.
%$$
%Let functions of the family $\varPhi$ be as in Theorem 3 and convex. 
Then $G({\Psi^*}) = GS({\varPhi^*})$.
\end{theorem}

The proof of Theorem 4 is essentially based on results of subsection 5.1 where some properties of the Young-Fenchel transform are considered.    

{\bf Remark 1}. Note that if functions $\varphi_m$ are defined on ${\mathbb R}^n$ by the formula 
$\varphi_m(x) = \Omega(2^m \Vert x \Vert)$ then the family $\varPhi$ satisfies the assumptions of Theorem~4.

\section{Auxiliary results}

In the proofs of the main results the following four lemmas will be useful. 

%Here we prove two assertions on properties of function  $\psi^*$ which will be used in the proofs of Theorems.

%\hspace {0,5cm} 

\begin{lemma} Let $g \in {\mathcal B}([0, \infty)^n)$. 
Then for each $M>0$ there exists a constant $A>0$ such that
$$
(g[e])^*(x) \le 
\sum \limits_{1 \le j \le n: x_j \ne 0} (x_j \ln\frac {x_j}{M} - x_j) + A, \ x \in [0, \infty)^n.
$$ 
\end{lemma}
%\overline \Pi

{\bf Proof}. For each $M>0$ we can find a number $A > 0$ 
such that for all $y = (y_1, \ldots , y_n) \in {\mathbb R}^n$ we have
$g[e](y) \ge \sum \limits_{j=1}^n M e^{y_j} - A$.
Hence, for $x = (x_1, \ldots , x_n) \in [0, \infty)^n$
$$
(g[e])^*(x)=\sup_{y \in {\mathbb R}^n}(\langle x, y \rangle - g[e](y))  \le 
\sup_{(y_1, \ldots , y_n) \in {\mathbb R}^n}\sum \limits_{j=1}^n (x_j y _j - M e^{y_j}) + A 
$$
$$
=\sum \limits_{1 \le j \le n: x_j \ne 0} \sup_{y_j \in {\mathbb R}}(x_jy_j -M e^{y_j}) + A = 
\sum \limits_{1 \le j \le n: x_j \ne 0} (x_j \ln\frac {x_j}{M} - x_j) + A. \ \square
$$ 
%From this the assertion of Lemma follows.

\begin{corollary} 
Let $g \in {\mathcal B}([0, \infty)^n)$. Then for each $b>0$ the series \ 
$\displaystyle \sum_{\vert j \vert \ge 0}^{\infty} \frac {e^{(g[e])^*(j)}} {b^{\vert j \vert} j!}$ and 
$\displaystyle \sum_{\vert j \vert \ge 0}^{\infty} \frac {e^{(g[e])^*(j)}} {b^{\vert j \vert} \vert j \vert!}$  are converging.
\end{corollary}

{\bf Remark 2}. Note that if $g \in {\mathcal B}([0, \infty)^n)$ then it is easy to see that 
$(g[e])^*(x)=+\infty$ for 
$x \notin [0, \infty)^n$, $(g[e])^*(x) > -\infty$ for $x \in {\mathbb R}^n$.
%This means, in particular, that  $dom \, {\psi_k^*} = [0, \infty)^n$ ($k \in {\mathbb N}$).
Also notice that $\displaystyle \lim_{x \to \infty, \atop x \in [0, \infty)^n} \frac {(g[e])^*(x)}{\Vert x \Vert}= + \infty$. Indeed, from the definition of $(g[e])^*$ it follows that 
$(g[e])^*(x) \ge \langle x, t \rangle - (g[e])(t)$
for all $x \in [0, \infty)^n$ and 
$t \in {\mathbb R}^n$. 
So if positive $M$ is arbitrary then from this inequality 
we get that $(g[e])^*(x) \ge M \Vert x \Vert - g[e](\frac {M x}{\Vert x \Vert})$ for $x \ne 0$. 
From this our assertion follows.

\begin{lemma} 
Let $u, v \in {\mathcal B}([0, \infty)^n)$ and for some $l > 0$ 
$$
2 u(x) \le v(x) + l, \ x \in [0, \infty)^n. 
$$

Then 
$$
(v[e])^*(x + y) \le (u[e])^*(x) + (u[e])^*( y) + l, \ x, y \in [0, \infty)^n.
$$
\end{lemma}
 
{\bf Proof}. 
Let $x, y \in [0, \infty)^n$. Then for each $t \in {\mathbb R}^n$ we have that
$$
(u[e])^*(x) + (u[e])^*( y) \ge \langle x+y, t \rangle - 2 u[e](t) \ge \langle x+y,  t \rangle  - v[e](t) - l.
$$
From this it follows that  
$$
(v[e])^*(x + y)  = \sup_{t \in {\mathbb R}^n} (\langle x+y, t \rangle - v[e](t)) \le (u[e])^*(x) + (u[e])^*( y)  + l. \  \square
$$

\begin{lemma} 
Let $u, v \in {\mathcal B}([0, \infty)^n)$ and there are constants $\sigma > 1$ and $\gamma > 0$ such that
$$
u(\sigma x) \le v(x) + \gamma, \ x \in [0, \infty)^n.
$$

Then for $x = (x_1, \ldots , x_n) \in [0, \infty)^n$ one has
$$
(u[e])^*(x) - (v[e])^*(x) \ge  \sum \limits_{j=1}^n x_j \ln \sigma - \gamma.
$$ 
\end{lemma}

{\bf Proof}. Note that by Lemma 1 and Remark 2 $(u[e])^*(x) < \infty$ and $(v[e])^*(x) < \infty$ for $x \in [0, \infty)^n$.
From the condition on $u$ and $v$ it follows that 
$$
u[e](t + \eta) \le v[e](t) + \gamma, \ t \in {\mathbb R}^n, 
$$
where $\eta = (\ln \sigma, \ldots , \ln \sigma)$.
Then for $x = (x_1, \ldots , x_n) \in [0, \infty)^n$
$$
(u[e])^*(x) - (v[e])^*(x) = \displaystyle \sup \limits_{t \in {\mathbb R}^n}(\langle x, t \rangle - u[e](t)) - 
\displaystyle \sup \limits_{t \in {\mathbb R}^n}(\langle x, t \rangle - v[e](t)) 
$$
$$
\ge 
\displaystyle \sup \limits_{t \in {\mathbb R}^n}(\langle x, t \rangle - u[e](t))  
- \displaystyle \sup \limits_{t \in {\mathbb R}^n}(\langle x, t \rangle - u[e](t+\eta)) - \gamma 
$$
$$
=
\displaystyle \sup \limits_{t \in {\mathbb R}^n}(\langle x, t \rangle - u[e](t)) - 
\displaystyle \sup \limits_{t \in {\mathbb R}^n}(\langle x, t+\eta \rangle - 
u[e](t+\eta)) + 
\langle x, \eta \rangle  - \gamma 
=  \sum \limits_{j=1}^n x_j \ln \sigma - \gamma. 
$$ 

\begin{lemma} Let $g = (g_1, \ldots , g_n)$ be a vector-function on ${\mathbb R}^n$ with convex components 
$g_j: {\mathbb R}^n \to [0, \infty)$ and a function $f: {\mathbb R}^n \to {\mathbb R}$ 
be such that $f_{|[0, \infty)^n}$ is convex and nondecreasing in each argument. Then $f \circ g$ is convex on ${\mathbb R}^n$.
\end{lemma}

{\bf Proof}. Let $x, y \in {\mathbb R}^n$, $\alpha \in [0, 1]$. 
%Then 
%using convexity of  
%for each  the function 
%$g_j$ on ${\mathbb R}^n$ 
Then for each $j = 1, \ldots , n$ we have that $0 \le g_j(\alpha x + (1 - \alpha) y) \le \alpha g_j(x) + (1 - \alpha) g_j(y)$. 
%Since $g_j(\alpha x + (1 - \alpha) y)
From this using monotonicity of $f$ on $[0, \infty)^n$ we get that 
$f(g(\alpha x + (1 - \alpha) y)) \le f(\alpha g(x) + (1 - \alpha) g(y))$. 
Now using the convexity of $f$ on $[0, \infty)^n$ we obtain the required relation 
$f(g(\alpha x + (1 - \alpha) y)) \le \alpha f(g(x)) + (1 - \alpha) f(g(y))$. $\square$

\section{Equivalent description of the space $E(\varPhi)$}

{\bf 3.1. Proof of Theorem 1}. 
Let $f \in E(\varPhi)$. 
Then $f \in E(\varphi_{\nu})$ for some $\nu \in {\mathbb N}$. 
Let $m \in {\mathbb Z_+}$, $\alpha = (\alpha_1, \ldots , \alpha_n) \in {\mathbb Z}_+^n$ and 
$x = (x_1, \ldots , x_n) \in {\mathbb R}^n$ be arbitrary. 
For $j =1, \ldots , n$ let $R_j$ be an arbitrary positive number. For $R = (R_1, \ldots , R_n)$ let 
$L_R(x)= \{\zeta = (\zeta_1, \ldots , \zeta_n) \in {\mathbb C}^n: \vert \zeta_j - x_j \vert = R_j, j=1, \ldots , n \}$.
Using Cauchy integral formula we have that
$$
(1+ \Vert x \Vert)^m (D^{\alpha}f)(x) =
\frac {\alpha! }{(2\pi i)^n} 
\displaystyle 
%\idotsint 
\int_{L_R(x)}
\frac 
{f(\zeta) (1+ \Vert x \Vert)^m \ d \zeta}
{(\zeta_1 - x_1)^{\alpha_1 +1} \cdots (\zeta_n - x_n)^{\alpha_n +1}} .
%d \zeta_1 \cdots d  \zeta_n ,
$$
From this we get that
$$
(1+ \Vert x \Vert)^m \vert (D^{\alpha}f)(x) \vert 
\le \frac {\alpha! }{(2\pi)^n} 
\displaystyle 
%\idotsint 
%\limits
\int_{L_R(x)}
\frac 
{(1+ \Vert x  - \zeta\Vert)^m (1+ \Vert \zeta \Vert)^m   \vert f(\zeta) \vert \, \vert d \zeta \vert}
{\vert \zeta_1 - x_1\vert^{\alpha_1 +1} \cdots \vert \zeta_n - x_n\vert^{\alpha_n +1}} 
$$
$$
\le \frac 
{\alpha! p_{\nu, m}(f) (1 + \Vert R \Vert)^m  e^{\varphi_{\nu}(R)}}{R^{\alpha}} .
$$
Using the condition $i_0$) on $\varPhi$ we obtain that
%for all $x \in {\mathbb R}^n$ and $R>0$
$$
(1+ \Vert x \Vert)^m \vert (D^{\alpha}f)(x) \vert \le 
e^{C_{\nu, m}} \alpha! p_{\nu, m}(f)
\frac 
{e^{\varphi_{\nu + 1}(R)}}{R^{\alpha}}.
$$
Hence, 
%for all $x \in {\mathbb R}^n$ 
$$
(1+ \Vert x \Vert)^m \vert (D^{\alpha}f)(x) \vert \le 
e^{C_{\nu, m}} \alpha! p_{\nu, m}(f)
\inf_{R \in (0, \infty)^n}
\frac 
{e^{\varphi_{\nu + 1}(R)}}{R^{\alpha}} 
$$
$$
= e^{C_{\nu, m}} \alpha! p_{\nu, m}(f)
\exp({-\sup_{r \in {\mathbb R}^n}
(\langle \alpha, r \rangle  - \psi_{\nu + 1} (r))}) = 
 e^{C_{\nu, m}} \alpha! p_{\nu, m}(f) e^{-\psi_{\nu + 1}^*(\alpha)}.
$$
From this it follows that for each $m \in {\mathbb Z}_+$   
\begin{equation}
\rho_{m, \nu + 1}(f_{|{\mathbb R}^n}) \le e^{C_{\nu, m}} p_{\nu, m}(f). 
\end{equation}
Therefore, 
$f_{|{\mathbb R}^n} \in {\mathcal E}(\psi_{\nu+1}^*)$.
Thus, $f_{|{\mathbb R}^n} \in {\mathcal E}(\Psi^*)$. Note that the inequality (1)
ensures the continuity of the embedding. $\square$

{\bf Proof of Theorem 2}. 
Let $f \in {\mathcal E}(\Psi^*)$. 
Then $f \in {\mathcal E}(\psi_{\nu}^*)$ for some $\nu \in {\mathbb N}$.
Hence, for each $m \in {\mathbb Z}_+$ we have that
%there exists a constant $d_m > 0$ such that for all $x \in {\mathbb R}^n$ and $\alpha \in {\mathbb Z}_+^n$ 
\begin{equation}
(1+ \Vert x \Vert)^m \vert (D^{\alpha}f)(x) \vert \le  \rho_{m, \nu} (f)
\alpha! e^{-\psi_{\nu}^*(\alpha)}, \ x \in {\mathbb R}^n, \alpha \in {\mathbb Z}_+^n.
\end{equation}
Since
$\displaystyle \lim_{x \to \infty, \atop x \in \Pi_n} \frac {\psi_{\nu}^*(x)}{\Vert x \Vert}= + \infty$ 
(see Remark 2)
then for each $\varepsilon >0$ there is a constant $c_{\varepsilon}>0$ such that  
for all $x \in {\mathbb R}^n$ and $\alpha \in {\mathbb Z}_+^n$ we have that
\begin{equation}
\vert (D^{\alpha}f)(x)\vert \le c_{\varepsilon} {\varepsilon}^{\vert \alpha \vert}\alpha!. 
\end{equation}
Hence, the sequence
$(\sum \limits_{\vert \alpha \vert \le k} \frac 
{(D^{\alpha}f)(0)}{\alpha!} x^{\alpha})_{k=1}^{\infty}$ 
converges to $f$ uniformly on compacts of ${\mathbb R}^n$.  
Also from (3) it follows that the series 
$
\displaystyle \sum_{\vert \alpha \vert \ge 0} \frac {(D^{\alpha}f)(0)}{\alpha!} z^{\alpha}
$
converges uniformly on compacts of ${\mathbb C}^{n}$ and, hence, 
its sum $F_f(z)$ is an entire function.
% in ${\mathbb C}^n$. 
Obviously, $F_f(x) = f(x), \ x \in {\mathbb R}^n$. The uniqueness of holomorphic continuation is obvious. 

Now we have to show that $F_f \in E(\varPhi)$. 
We will estimate a growth  of $F_f$ using the inequality (2) and the Taylor series expansion 
of $F_f(z)$ ($z = x+iy, x, y \in {\mathbb R}^n$) with respect to a point $x$: 
$$
F_f(z) = \displaystyle \sum_{\vert \alpha \vert \ge 0} \frac {(D^{\alpha}f)(x)}{\alpha!} (iy)^{\alpha}.
$$
%where $z = x+y$ $ (y \in {\mathbb R}^n)$  
Let $m \in {\mathbb Z}_+$ be arbitrary. Then
$$
(1 + \Vert z \Vert)^m \vert F_f(z) \vert \le \sum_{\vert \alpha \vert \ge 0} \frac {1}{\alpha!}
(1 + \Vert x \Vert)^m (1 + \Vert y \Vert)^m
\prod \limits_{j=1}^n(1 + \vert y_j \vert)^{\alpha_j} \vert (D^{\alpha}f)(x)\vert  
$$
$$
\le \rho_{m, \nu} (f) (1 + \Vert y \Vert)^m
\sum_{\vert \alpha \vert \ge 0}
e^{-\psi_{\nu}^*(\alpha)}\prod \limits_{j=1}^n(1 + \vert y_j \vert)^{\alpha_j}  
$$
$$
\le
\rho_{m, \nu} (f) 
(1 + \Vert y \Vert)^m
\sum_{\vert \alpha \vert \ge 0} 
\frac {\prod \limits_{j=1}^n(1 + \vert y_j \vert)^{\alpha_j}}
{e^{\psi^*_{\nu+1}(\alpha)}}e^{\psi^*_{\nu+1}(\alpha)-\psi^*_{\nu}(\alpha)}.
$$
Since $\varPhi$ satisfies the condition $i_1)$ then by Lemma 3 
$$
\psi_{\nu}^*(x) - \psi_{\nu + 1}^*(x) \ge \delta_{\nu} \sum \limits_{j=1}^n x_j - \gamma_{\nu}, \ 
x = (x_1, \ldots , x_n) \in  [0, \infty)^n,
$$
where $\delta_{\nu} = \ln \sigma_{\nu}$.
Using this inequality and denoting 
the sum of the series 
$\displaystyle\sum_{\vert \alpha \vert \ge 0} e^{\psi^*_{\nu+1}(\alpha)-\psi^*_{\nu}(\alpha)}$
by $B_{\nu}$ we have that
$$
(1 + \Vert z \Vert)^m \vert F_f(z) \vert \le B_{\nu} \rho_{m, \nu} (f) (1 + \Vert y \Vert)^m 
\sup \limits_{\vert \alpha \vert \ge 0}
\frac {\prod \limits_{j=1}^n(1 + \vert y_j \vert)^{\alpha_j}}{e^{\psi^*_{\nu+1}(\alpha)}}
$$
$$
\le B_{\nu} \rho_{m, \nu} (f) 
(1 + \Vert y \Vert)^m 
e^{\sup \limits_{t = (t_1, \ldots , t_n) \in [0, \infty)^n}
(t_1 \ln (1 + \vert y_1 \vert) + \cdots + t_n \ln (1 + \vert y_n \vert) -  \psi^*_{\nu+1}(t))}.
$$ 
Taking into account Remark 2 we obtain that
\begin{equation}
(1 + \Vert z \Vert)^m \vert F_f(z) \vert 
%\le c_m B_{\nu} (1 + \vert y \vert)^m e^{\sup \limits_{t \ge 0}(t \ln (1 + \vert y \vert) -  \psi^*_{\nu+1}(t))} = 
\le B_{\nu} \rho_{m, \nu} (f) 
e^{(\psi^*_{\nu+1})^*(\ln (1 + \vert y_1 \vert), \ldots , \ln (1 + \vert y_n \vert)) + m \ln (1 + \Vert y \Vert))}.
\end{equation}

Note that for each $k \in {\mathbb N}$ and $A > 0$ from condition $i_0)$ we can find a constant $C_1(k, A) > 0$ such that for all 
$x = (x_1, \ldots , x_n) \in {\mathbb R}^n$
$$
\psi_k(x) + A \displaystyle \sum_{j=1}^n x_j \le \psi_{k+1}(x) + C_1(k, A).
$$
So for all $\xi \in [0, \infty)^n$ and $A>0$
$$
\psi_k^*(\xi) = 
\displaystyle \sup \limits_{x \in {\mathbb R}^n}
(\langle \xi,  x \rangle - \psi_k(x))  \ge 
\displaystyle \sup \limits_{x \in {\mathbb R}^n}(\langle \xi,  x \rangle - \psi_{k + 1}(x) + 
A \displaystyle \sum_{j=1}^n x_j) - C_1(k, A).
$$
Put $\Lambda:=(A, \ldots , A) \in {\mathbb R}^n$. Then from the previous estimate we have that 
$$
\psi_k^*(\xi) \ge 
\displaystyle \sup \limits_{x \in {\mathbb R}^n}(\langle \xi + \Lambda, x \rangle - \psi_{k + 1}(x)) - C_1(k, A) = 
\psi_{k + 1}^*(\xi + \Lambda) - C_1(k, A).
$$
Further, for all $x \in [0, \infty)^n$ and $A > 0$
$$
(\psi_k^*)^*(x) =  \displaystyle \sup \limits_{\xi \in {\mathbb R}^n}
(\langle x, \xi \rangle - \psi_k^*(\xi)) = \displaystyle \sup \limits_{\xi \in [0, \infty)^n}
(\langle x, \xi \rangle - \psi_k^*(\xi)) 
$$
$$
\le
\displaystyle \sup \limits_{\xi \in [0, \infty)^n}
(\langle x, \xi \rangle  - \psi_{k + 1}^*(\xi + \Lambda)) + C_1(k, A)
$$
$$
= \displaystyle \sup \limits_{\xi \in [0, \infty)^n}
(\langle x, \xi + \Lambda \rangle  - \psi_{k + 1}^*(\xi + \Lambda)) - A \displaystyle \sum_{j=1}^n x_j + C_1(k, A) 
$$
$$
\le \displaystyle \sup \limits_{\xi \in [0, \infty)^n}
(\langle x, \xi \rangle  - \psi_{k + 1}^*(\xi)) - A \displaystyle \sum_{j=1}^n x_j + C_1(k, A) 
$$
$$=
(\psi_{k + 1}^*)^*(x) - A \displaystyle \sum_{j=1}^n x_j + C_1(k, A).
$$
Thus, for each $k \in {\mathbb N}$ and $A > 0$ we have
%there exists a constant $C(m, A) > 0$ such that 
\begin{equation}
(\psi_k^*)^*(x) + A \displaystyle \sum_{j=1}^n x_j \le (\psi_{k+1}^*)^*(x) + C_1(k, A), \ x \in [0, \infty)^n.
\end{equation}
Now using the inequality  (5) for $A=m$ and $k=\nu +1$, we obtain from the estimate (4) that
\begin{equation}
(1 + \Vert z \Vert)^m \vert F_f(z) \vert 
\le B_{\nu} \rho_{m, \nu} (f)  
e^{C_1(\nu + 1, m)} 
e^{(\psi_{\nu + 2}^*)^* (\ln (1 + \vert y_1 \vert), \ldots , \ln (1 + \vert y_n \vert))}.
\end{equation}
Since $(\psi_{\nu + 2}^*)^*(t) \le \psi_{\nu + 2}(t)$ for $t \in [0, \infty)^n$ then from (6) we get that 
$$
(1 + \Vert z \Vert)^m \vert F_f(z) \vert 
\le B_{\nu} \rho_{m, \nu} (f)  e^{C_1(\nu + 1, m)} e^{\psi_{\nu + 2}(\ln (1 + \vert y_1 \vert), \ldots , \ln (1 + \vert y_n \vert))}.
$$
In other words,
$$
(1 + \Vert z \Vert)^m \vert F_f(z) \vert 
\le B_{\nu}
\rho_{m, \nu} (f)  e^{C_1(\nu + 1, m)} e^{\varphi_{\nu + 2} (1 + \vert y_1 \vert, \ldots , 1 + \vert y_n \vert)}.
$$
Using the condition $i_2)$ on $\varPhi$ it is possible to find a constant $K_{\nu, m} > 0$ 
such that for all $z \in {\mathbb C}^n$
\begin{equation}
(1 + \Vert z \Vert)^m \vert F_f(z) \vert \le K_{\nu, m}
\rho_{m, \nu} (f)  e^{\varphi_{\nu+3} (\vert Im z_1 \vert, \ldots , \vert Im z_n \vert)}.
\end{equation} 
Thus, for each $m \in {\mathbb Z}_+$ we have that $p_{\nu + 3, m} (F_f) \le K_{\nu, m} \rho_{m, \nu} (f)$. 
Hence, 
$F_f \in E(\varphi_{\nu + 3})$. Thus, $F_f \in E(\varPhi)$. Also note that the last inequality 
ensures the continuity of the embedding. $\square$

{\bf 3.2.  Another structure of $E(\varPhi)$}.
For each $\nu \in {\mathbb N}$ and $m \in {\mathbb Z}_+$ consider the normed space
$$
{\mathcal H}_m(\varphi_{\nu})= \{f \in H({\mathbb C}^n): 
\sigma_{{\nu}, m}(f) = 
\sup_{z \in {\mathbb C}^n} 
\frac 
{\vert f(z)\vert (1 + \Vert z \Vert)^m}
{e^{(\psi_{\nu}^*)^*(\ln (1 + \vert Im z_1 \vert), \ldots , \ln (1 + \vert Im z_n \vert ))}} < \infty \}.
$$
Let ${\mathcal H}(\varphi_{\nu})= \bigcap \limits_{m=0}^{\infty} {\mathcal H}_m(\varphi_{\nu})$.
%, $E(\varPhi)= \bigcup \limits_{\nu=1}^{\infty} E(\varphi_{\nu})$. 
%Since $p_{\nu, k}(f) \le p_{\nu, k+1}(f)$ 
%for $f \in E_{k+1}(\varphi_{\nu})$ 
%then 
Obviously, ${\mathcal H}_{m+1}(\varphi_{\nu})$ is continuously embedded in ${\mathcal H}_m(\varphi_{\nu})$. 
Endow ${\mathcal H}(\varphi_{\nu})$ with a projective limit topology of spaces ${\mathcal H}_m(\varphi_{\nu})$. 
Note that if $f \in {\mathcal H}(\varphi_{\nu})$ then using the inequality (5) we have that
$\sigma_{\nu+1, m}(f) \le e^{C_1(\nu, 1)} \sigma_{\nu, m}(f)$ for each $m \in {\mathbb Z}_+$.
Thus, ${\mathcal H}(\varphi_{\nu})$ is continuously embedded in ${\mathcal H}(\varphi_{\nu + 1})$ for each 
$\nu \in {\mathbb N}$. 
%Let ${\mathcal H}(\varPhi)= \bigcup \limits_{\nu=1}^{\infty} {\mathcal H}(\varphi_{\nu})$. With the usual operations %of addition and multiplication by complex numbers 
%${\mathcal H}(\varPhi)$
%and $E(\varPhi)$ are 
%is a linear space. 
Supply ${\mathcal H}(\varPhi)= \bigcup \limits_{\nu=1}^{\infty} {\mathcal H}(\varphi_{\nu})$ with the topology of the inductive limit of spaces ${\mathcal H}(\varphi_{\nu})$.

\begin{proposition}
Let all the functions of the family $\varPhi$ satisfy the condition $i_2)$ of Theorem 2 and every function $\psi_{\nu}$ be convex on ${\mathbb R}^n$ ($\nu \in {\mathbb N}$). 
Then $E(\varPhi) = {\mathcal H}(\varPhi)$. 
\end{proposition}

{\bf Proof}. By assumption each function $\psi_{\nu}$ is convex and continuous on ${\mathbb R}^n$. Since
the Young-Fenchel conjugation is involutive (see \cite {R}, Theorem 12.2) it follows that 
$(\psi_{\nu}^*)^* = \psi_{\nu}$. 
Thus, for each $\nu \in {\mathbb N}$ and $t = (t_1, \ldots , t_n) \in [0, \infty)^n$ we have that
$$
(\psi_{\nu}^*)^*(\ln (1 + t_1), \ldots , \ln (1 + t_n))=\psi_{\nu}(\ln (1 + t_1), \ldots , \ln (1 + t_n))
$$
$$
=
\varphi_{\nu}(1 + t_1, \ldots , 1+t_n). 
$$
From this, and taking into account that functions of the family $\Phi$ are nondecreasing 
in each variable in $[0, \infty)^n$, we get that 
%for all $t = (t_1, \ldots , t_n) \in \overline [0, \infty)^n$ 
$$
(\psi_{\nu}^*)^*(\ln (1 + t_1), \ldots , \ln (1 + t_n)) \ge \varphi_{\nu}(t), \ 
t = (t_1, \ldots , t_n) \in [0, \infty)^n.
$$
On the other hand using the condition $i_2)$ of Theorem 2 we have 
\begin{equation}
(\psi_{\nu}^*)^*(\ln (1 + t_1), \ldots , \ln (1 + t_n)) \le \varphi_{\nu + 1}(t)  + K_{\nu}, \ t \in [0, \infty)^n. 
\end{equation}
From these inequalities the assertion follows. $\square$

Using Theorems 1 and 2 we can prove the following

\begin{proposition}
Let the family $\varPhi$ satisfies the conditions of Theorem 2. 
Then $E(\varPhi) = {\mathcal H}(\varPhi)$.
\end{proposition}

{\bf Proof}. Thanks to the condition $i_2)$ the inequality (8) holds. Using it we have that 
for each $m \in {\mathbb Z}_+$ 
$$
p_{\nu + 1, m}(f) \le e^{K_{\nu}} \sigma_{{\nu}, m}(f), \ f \in {\mathcal H}(\varphi_{\nu}).
$$
Hence, the identity embedding $I: {\mathcal H}(\varPhi) \to E(\varPhi)$ is continuous. 
%Let $f \in {\mathcal H}(\varPhi)$. Hence, 
%$f \in {\mathcal H}(\varphi_{\nu})$ for some $\nu \in {\mathbb N}$. 
%Then for each $k \in {\mathbb Z}_+$ we have that 
%$$
%\vert f(z) \vert \le {\cal N}_{\varphi_{\nu}, k}(f) 
%\frac {e^{((\varphi_{\nu}[e])^*)^*(\ln (1 + \vert Im z_1 \vert), \ldots , \ln (1 + \vert Im z_n \vert ))}} 
%{(1 + \Vert z \Vert)^k} \le 
%$$
%$$
%\le {\cal N}_{\varphi_{\nu}, k}(f)
%\frac {e^{\varphi_{\nu}[e](\ln (1 + \vert Im z_1 \vert), \ldots , \ln (1 + \vert Im z_n \vert ))}}
%{(1 + \Vert z \Vert)^k} = {\cal N}_{\varphi_{\nu}, k}(f)
%\frac {e^{\varphi_{\nu}(1 + \vert Im z_1 \vert, \ldots , 1 + \vert Im z_n \vert)}}
%{(1 + \Vert z \Vert)^k} .
%$$
%From this using the condition $i_2)$ we have that for each $k \in {\mathbb Z}_+$ 
%$$
%\vert f(z) \vert \le e^{K_{\nu}} {\cal N}_{\varphi_{\nu}, k}(f) 
%\frac {e^{\varphi_{\nu}(\vert Im z_1 \vert, \ldots , \vert Im z_n \vert)}}
%{(1 + \Vert z \Vert)^k} .
%$$
%Therefore, if $f \in {\mathcal H}(\varphi_{\nu})$ then for each $k \in {\mathbb Z}_+$
%$$
%p_{\nu + 1, k}(f) \le e^{K_{\nu}} {\cal N}_{\varphi_{\nu}, k}(f), \ f \in {\mathcal H}(\varphi_{\nu}).
%$$
%Hence, $f \in E(\varphi_{\nu + 1})$. Thus, $f \in E(\varPhi)$. 
%Also note that from the last inequality it follows that the embedding $I: {\mathcal H}(\varPhi) \to E(\varPhi)$ is %continuous. 

The mapping $I$ is surjective too. Indeed, if $f \in E(\varPhi)$ then $f \in E(\varphi_{\nu})$ for some $\nu \in {\mathbb N}$. 
Let $m \in {\mathbb Z}_+$ be arbitrary. 
By the inequality (1)  we have that  
$
\rho_{m, \nu + 1}(f_{|{\mathbb R}^n}) \le e^{C_{\nu, m}} p_{\nu, m}(f).
$
%Also taking into account 
From this and the inequality (6) (with $\nu$ replaced by $\nu + 1$) we obtain that
$$
\sigma_{{\nu + 3}, m}(f) \le A_{\nu, m} p_{\nu, m}(f),
$$
where $A_{\nu, m}$ is some positive constant. Hence, $f \in {\mathcal H}(\varphi_{\nu + 3})$. 
So, $f \in {\mathcal H}(\varPhi)$. Moreover, the last estimate shows that the inverse mapping $I^{-1}$ is continuous.
Hence, the equality $E(\varPhi) = {\mathcal H}(\varPhi)$ is topological too. $\square$

\section{Fourier transform of $E(\varPhi)$}

Recall two notations which will be used in the proof of Theorem 3. Namely, for 
$\alpha = ({\alpha}_1, \ldots , {\alpha}_n)$ and 
$\beta = ({\beta}_1, \ldots , {\beta}_n) \in {\mathbb Z_+^n}$ the notation
$\alpha \le \beta $ indicates that
${\alpha}_j \le {\beta}_j$ ($j = 1, 2, \ldots , n$) 
and in such case $\binom {\beta}{\alpha}:= \prod \limits_{j=1}^{n} \binom {\beta_j}{\alpha_j}$,
%we denote by
%$C_{\beta}^{\alpha}:= \prod \limits_{j=1}^{n} C_{\beta_j}^{\alpha_j}$, 
where $\binom {\beta_j}{\alpha_j}$
%$C_{\beta_j}^{\alpha_j}$ 
are the binomial coefficients. 

{\bf Proof of Theorem 3}. Let $\nu \in {\mathbb N}$ and $f \in E(\varphi_{\nu})$.
Let $\alpha$, $\beta \in {\mathbb Z}_+^n$, $x, \eta \in {\mathbb R}^n$. Then
$$
x^{\beta} ({D^{\alpha} \hat f)(x) = x^{\beta} \int_{{\mathbb R}^n}} f(\zeta) 
(-i \zeta)^{\alpha} 
e^{-i \langle x, \zeta \rangle} \ d \xi, \  \zeta = \xi + i\eta.
$$
From this equality we have that
$$
\vert x^{\beta} (D^{\alpha} \hat f)(x)\vert \le 
\int_{{\mathbb R}^n}
\vert f(\zeta) \vert
\Vert \zeta \Vert^{\vert \alpha \vert} 
e^{\langle x, \eta \rangle} \prod \limits_{j=1}^n \vert x_j \vert^{\beta_j} \ d \xi  
$$
$$
\le 
\int_{{\mathbb R}^n} 
\vert f(\zeta) \vert
(1 + \Vert \zeta \Vert)^{n + \vert \alpha \vert + 1} 
e^{\langle x, \eta \rangle} \prod \limits_{j=1}^n \vert x_j \vert^{\beta_j} \ 
\frac {d \xi}{(1 + \Vert \xi \Vert)^{n+1}} \ .
$$
Hence, 
\begin{equation}
\vert x^{\beta} (D^{\alpha} \hat f)(x)\vert
\le 
s_n p_{\nu, n + \vert \alpha \vert + 1}(f) e^{\varphi_{\nu}(\eta)}
e^{\langle x, \eta \rangle} \prod \limits_{j=1}^n \vert x_j \vert^{\beta_j} .
\end{equation}
If $\vert \beta \vert = 0$ then (putting $\eta = 0$ in (9)) we have that
\begin{equation}
\vert (D^{\alpha}\hat f)(x) \vert 
\le  s_n e^{\varphi_{\nu} (0)} p_{\nu, n + \vert \alpha \vert + 1}(f).
\end{equation} 
Now let consider the case when $\vert \beta \vert > 0$.
For $x =(x_1, \ldots , x_n) \in {\mathbb R}^n$ let $\theta(x)$ be a point in ${\mathbb R}^n$ with coordinates $\theta_j$ defined as follows: $\theta_j=\frac {x_j} {\vert x_j\vert}$ if $x_j \ne 0$ and $\theta_j=0$ if $x_j=0$ ($j=1, \ldots , n)$.
Let $t = (t_1, \ldots, t_n) \in {\mathbb R}^n$ have strictly positive coordinates.
Put
$\eta = -(\theta_1 t_1, \ldots , \theta_n t_n)$.
%, $\tilde x =(\vert x_1 \vert, \ldots , \vert x_n \vert)$. 
Then  from (9) we get that
$$
\vert x^{\beta} 
(D^{\alpha}\hat f)(x)\vert \le 
s_n p_{\nu, n + \vert \alpha \vert + 1}(f) 
e^{\varphi_{\nu}(t)} \prod \limits_{j \in \{1, \ldots , n\}: \beta_j \ne 0} 
\frac {\vert x_j \vert^{\beta_j}}{e^{t_j \vert x_j \vert}}  
$$
$$
\le
s_n p_{\nu, n + \vert \alpha \vert + 1}(f) 
e^{\varphi_{\nu}(t)} \prod \limits_{j \in \{1, \ldots , n\}: \beta_j \ne 0} 
\sup \limits_{r_j > 0}
\frac {r_j^{\beta_j}}{e^{t_j r_j}}
$$
$$
=
s_n p_{\nu, n + \vert \alpha \vert + 1}(f)
\exp (\varphi_{\nu}(t) + \sum \limits_{1 \le j \le n: \beta_j \ne 0} \sup \limits_{r_j > 0}(- t_j r_j + \beta_j \ln r_j))
$$
$$
=
s_n p_{\nu, n + \vert \alpha \vert + 1}(f)\exp (\varphi_{\nu}(t) + \sum \limits_{1 \le j \le n: \beta_j \ne 0} (\beta_j \ln \beta_j - \beta_j) - \sum \limits_{j=1}^n \beta_j \ln t_j). 
$$
From this we have that
$$
\vert x^{\beta} 
(D^{\alpha}\hat f)(x)\vert \le 
s_n p_{\nu, n + \vert \alpha \vert + 1}(f) 
e^{\sum \limits_{1 \le j \le n: \atop \beta_j \ne 0} \beta_j \ln \frac {\beta_j}{e}  + \inf \limits_{t = (t_1, \ldots , t_n) \in (0, \infty)^n} 
(\varphi_{\nu}(t) - \sum \limits_{j=1}^n \beta_j \ln t_j)}. 
$$
Note that for each  $\mu = (\mu_1, \ldots , \mu_n) \in {\mathbb Z}_+^n$
$$
\inf \limits_{t = (t_1, \ldots , t_n) \in (0, \infty)^n} (- \mu_1 \ln t_1 - \cdots - \mu_n \ln t_n + 
\varphi_{\nu}(t))
$$
$$
= - \sup_{u \in{\mathbb R}^n}(\langle \mu, u \rangle - \psi_{\nu}(u)) = - \psi_{\nu}^*(\mu),
$$
From this and the previous estimate it follows that
$$
\vert x^{\beta} 
(D^{\alpha} \hat f)(x)\vert \le s_n
p_{\nu, n + \vert \alpha \vert + 1}(f)  
e^{\sum \limits_{1 \le j \le n: \beta_j \ne 0} \beta_j \ln \frac {\beta_j}{e}} 
e^{- \psi_{\nu}^*(\beta)}.
$$
Now from this and (10) we have that for $\alpha, \beta \in {\mathbb Z}_+^n$, 
$x \in {\mathbb R}^n$
$$
\vert x^{\beta} 
(D^{\alpha} \hat f)(x)\vert \le s_n
p_{\nu, n + \vert \alpha \vert + 1}(f) \beta! e^{- \psi_{\nu}^*(\beta)}.
$$
From this it follows that for each $m \in {\mathbb Z}_+$
$$
\max_{\vert \alpha \vert \le m} \sup_{x \in {\mathbb R}^n, \beta \in {\mathbb Z}^n_+}
\frac {\vert x^{\beta} 
(D^{\alpha} \hat f)(x)\vert}{\beta! e^{- \psi_{\nu}^*(\beta)}} \le
s_n p_{\nu, n + m + 1}(f), \ f \in E(\varphi_{\nu}).
$$
In other words, for each $m \in {\mathbb Z}_+$
$$
\Vert \hat f \Vert_{m, \psi_{\nu}^*} \le s_n p_{\nu, n + m + 1}(f), \ f \in E(\varphi_{\nu}).
$$
From this inequality it follows that the linear mapping ${\cal F}: f \in E(\varPhi)  \to \hat f$ acts from  
$E(\varPhi)$ to $G(\Psi^*)$ and is continuous. 

Let us show that ${\mathcal F}$ is surjective. 
Take $g \in G(\Psi^*)$. 
Then $g \in G(\psi_{\nu}^*)$ for some $\nu \in {\mathbb N}$. Let $m \in {\mathbb Z}_+$ be 
arbitrary. Then for all $\alpha \in {\mathbb Z_+^n}$ with $\vert \alpha \vert \le m$, $\beta \in {\mathbb Z_+^n}$, 
$x \in {\mathbb R}^n$
we have that
$$
\vert x^{\beta}(D^{\alpha}g)(x) \vert \le 
\Vert g \Vert_{m, \psi^*_{\nu}} \beta! e^{-\psi_{\nu}^*(\beta)}.
$$
Using this inequality and the equality
$$
\vert (D^{\alpha}g)(x) \vert\prod \limits_{k=1}^n (1 + \vert x_k \vert)^{\beta_k}  = \vert (D^{\alpha}g)(x) \vert
\prod \limits_{k=1}^n \limits \sum \limits_{j_k=0}^{\beta_k} \binom {\beta_k}{j_k} \vert x_k \vert^{j_k}
$$
$$
=
\sum \limits_{j \in {\mathbb Z^n_+}: \atop (0, \ldots , 0) \le j \le \beta} \binom {\beta} {j} \vert x^j (D^{\alpha}g)(x) \vert 
$$
%Note that for all $\alpha \in {\mathbb Z_+^n}$, $\beta = (\beta_1, \ldots , \beta_n) \in {\mathbb Z_+^n}$, $x = %(x_1, \ldots , x_n) \in {\mathbb R}^n$ 
%\le 
%\Vert g \Vert_{m, \psi^*_{\nu}} \beta! e^{-\psi_{\nu}^*(\beta)}.
%$$
we have that 
%for all $\alpha \in {\mathbb Z_+^n}$ with $\vert \alpha \vert \le m$, $\beta \in {\mathbb Z_+^n}$
%$x \in {\mathbb R}^n$
\begin{equation}
\vert (D^{\alpha}g)(x) \vert\prod \limits_{k=1}^n (1 + \vert x_k \vert)^{\beta_k} 
\le \Vert g \Vert_{m, \psi^*_{\nu}}
\sum \limits_{j \in {\mathbb Z^n_+}: \atop (0, \ldots , 0) \le j \le \beta} \binom {\beta}{j} 
j! e^{-\psi_{\nu}^*(j)}.
\end{equation}
Now note that since the family $\varPhi$ satisfies  the condition $i_4)$ then 
with the help of Lemma 2 we have that for each $k \in {\mathbb N}$ 
%the following inequality holds
\begin{equation}
\psi_{k+1}^*(x+y) \le \psi_k^*(x) + \psi_k^*(y) + l_k, \ x, y \in [0, \infty)^n.  
\end{equation}
Using this inequality we have from (11) that 
$$
\vert (D^{\alpha}g)(x) \vert\prod \limits_{k=1}^n (1 + \vert x_k \vert)^{\beta_k} \le 
\Vert g \Vert_{m, \psi^*_{\nu}} e^{-\psi_{\nu + 1}^*(\beta) + l_{\nu}}
\sum \limits_{j \in {\mathbb Z^n_+}: \atop (0, \ldots , 0) \le j \le \beta} \binom {\beta}{j} 
j!  e^{\psi_{\nu}^*(\beta - j)}.
$$
From this we obtain that
$$
\vert (D^{\alpha}g)(x) \vert\prod \limits_{k=1}^n (1 + \vert x_k \vert)^{\beta_k} \le 
\Vert g \Vert_{m, \psi^*_{\nu}} \beta! e^{-\psi_{\nu + 1}^*(\beta) + l_{\nu}}
\sum \limits_{j \in {\mathbb Z^n_+}: \atop (0, \ldots , 0) \le j \le \beta} 
\frac {e^{\psi_{\nu}^*(\beta - j)}}{(\beta - j)!} .
$$ 
Recall that by the Corollary 1 the series 
$\sum \limits_{j \in {\mathbb Z^n_+}} 
\frac {e^{\psi_{\nu}^*(j)}}{j!}$ is converging. 
From this and the previous inequality it follows that 
for all $\alpha \in {\mathbb Z_+^n}$ with $\vert \alpha \vert \le m$, $\beta = (\beta_1, \ldots , \beta_n) \in {\mathbb Z_+^n}$, 
$x = (x_1, \ldots , x_n) \in {\mathbb R}^n$ 
\begin{equation}
\vert (D^{\alpha}g)(x) \vert \prod \limits_{k=1}^n (1 + \vert x_k \vert)^{\beta_k} \le c_1
\Vert g \Vert_{m, \psi^*_{\nu}} \beta! e^{-\psi_{\nu + 1}^*(\beta)} ,
\end{equation} 
where $c_1 = e^{l_{\nu}}\sum \limits_{j \in {\mathbb Z^n_+}} 
\frac {e^{\psi_{\nu}^*(j)}}{j!} $.
Now let
$$
f(\xi) = 
\frac {1}{(2 \pi)^n} 
\int_{{\mathbb R}^n} g(x) e^{i \langle x, \xi \rangle} \ dx, \ 
\xi \in {\mathbb R}^n.
$$
For all $\alpha = (\alpha_1, \ldots , \alpha_n), \beta = (\beta_1, \ldots , \beta_n) \in {\mathbb Z_+^n}, \xi \in {\mathbb R}^n$ we have that
$$
(i\xi)^{\beta}(D^{\alpha}f)(\xi) = \frac {(-1)^{\vert \beta \vert} }{(2 \pi)^n} 
\int_{{\mathbb R}^n} D^{\beta}(g(x) (ix)^{\alpha}) e^{i \langle x, \xi \rangle} \ dx .
$$
For $s=1, \ldots , n$ put $\gamma_s =\min (\beta_s, \alpha_s)$ and take $\gamma = (\gamma_1, \ldots , \gamma_n)$.
Then 
$$
(i\xi)^{\beta}(D^{\alpha}f)(\xi) = \frac {(-1)^{\vert \beta \vert} }{(2 \pi)^n} 
\int_{{\mathbb R}^n}
\displaystyle \sum \limits_{j \in {\mathbb Z_+^n}: j \le \gamma} \binom {\beta}{j}
 (D^{\beta - j} g)(x) (D^j (ix)^{\alpha}) 
e^{i \langle x, \xi \rangle} \ dx .
$$
From this we have that
$$
\vert \xi^{\beta}(D^{\alpha}f)(\xi) \vert \le \frac {1}{(2 \pi)^n}
\displaystyle \sum \limits_{j \in {\mathbb Z_+^n}: j \le \gamma}
\binom {\beta}{j}
\int_{{\mathbb R}^n} \vert(D^{\beta - j} g)(x)\vert 
\frac {\alpha!}{(\alpha - j)!} 
\vert x^{\alpha - j} \vert  \ dx
$$ 
$$
\le \frac {1}{(2 \pi)^n}
\displaystyle \sum \limits_{j \in {\mathbb Z_+^n}: j \le \gamma}
\frac {\binom {\beta}{j} \alpha!}{(\alpha - j)!}
\int \limits_{{\mathbb R}^n} \vert(D^{\beta - j} g)(x)\vert 
\prod \limits_{k=1}^n (1 + \vert x_k \vert)^{\alpha_k  - j_k + 2} 
\frac {dx}{\prod \limits_{k=1}^n (1 + \vert x_k \vert)^2} .
$$
Using the inequality (13) and denoting the element $(2, \ldots , 2) \in {\mathbb R}^n$ by $\omega$ we have that
$$
\vert \xi^{\beta}(D^{\alpha}f)(\xi) \vert \le \frac 
{c_1}{2^n} \Vert g \Vert_{\vert \beta \vert, \psi^*_{\nu}}
\displaystyle \sum \limits_{j \in {\mathbb Z_+^n}: j \le \gamma}
\frac {\binom {\beta}{j} \alpha!}{(\alpha - j)!} (\alpha - j + \omega)! e^{-\psi_{\nu + 1}^*(\alpha - j + \omega)}.
$$
Note that using the condition $i_4)$ it is easy to verify that 
$\varkappa_{\nu}: = \sup \limits_{x \in {\mathbb R}^n} (\psi_{\nu+2}^*(x) - \psi_{\nu+1}^*(x + \omega)) < \infty$. 
From this and the previous inequality it follows that
$$
\vert \xi^{\beta}(D^{\alpha}f)(\xi) \vert \le {c_2} \Vert g \Vert_{\vert \beta \vert, \psi^*_{\nu}}(\alpha + \omega)!
\displaystyle \sum \limits_{j \in {\mathbb Z_+^n}: j \le \gamma}
\binom {\beta}{j}  e^{-\psi_{\nu + 2}^*(\alpha - j)},
$$
where $c_2 = \frac {c_1 e^{\varkappa_{\nu}}}{2^n}$. 
Using the inequality (12) we get
$$
\vert \xi^{\beta}(D^{\alpha}f)(\xi) \vert \le 
c_2 e^{l_{\nu + 2}} \Vert g \Vert_{\vert \beta \vert, \psi^*_{\nu}}(\alpha + \omega)! 
e^{-\psi_{\nu + 3}^*(\alpha)} 
\displaystyle \sum \limits_{j \in {\mathbb Z_+^n}: j \le \gamma}
\binom {\beta}{j}  e^{\psi_{\nu + 2}^*(j)}.
$$
From this we obtain that
$$
\vert \xi^{\beta}(D^{\alpha}f)(\xi) \vert \le 
c_2 e^{l_{\nu + 2}} \Vert g \Vert_{\vert \beta \vert, \psi^*_{\nu}}(\alpha + \omega)! 
e^{-\psi_{\nu + 3}^*(\alpha)}  \beta!
\displaystyle \sum \limits_{j \in {\mathbb Z_+^n}: j \le \gamma}
\frac 
{e^{\psi_{\nu + 2}^*(j)}}{j!}.
$$
Take into account that for all $m_1, m_2 \in {\mathbb Z_+}$ 
$$
(m_1 + m_2)! \le e^{m_1 + m_2} m_1! m_2! 
$$
(this inequality easily follows from the inequality $(m_1 + m_2)^{m_2} \le m_2! e^{m_1 + m_2}$).
Using this inequality we get from the preceding inequality that
$$
\vert \xi^{\beta}(D^{\alpha}f)(\xi) \vert \le c_3 e^{\vert \alpha \vert} \beta! 
\Vert g \Vert_{\vert \beta \vert, \psi^*_{\nu}} \alpha!
e^{-\psi_{\nu + 3}^*(\alpha)},
$$
where 
$c_3 = c_2 2^n e^{2n + l_{\nu +2}}  \displaystyle \sum \limits_{j \in {\mathbb Z_+^n}}
\frac {e^{\psi_{\nu + 2}^*(j)}}{j!}$.
%Since $\varPhi$ satisfies 
From this using the inequality 
%we have that for each $k \in {\mathbb N}$
$$
\psi_k^*(x) - \psi_{k + 1}^*(x) \ge \sum \limits_{j=1}^n x_j \ln 2 - a_k, \ x \in [0, \infty)^n, k \in {\mathbb N},
$$
(that holds in view of the condition $i_3)$ and Lemma 3) 
%Using this inequality 
we obtain that
$$
\vert \xi^{\beta}(D^{\alpha}f)(\xi) \vert \le c_4 \beta!  
 \Vert g \Vert_{\vert \beta \vert, \psi^*_{\nu}} \alpha! e^{-\psi_{\nu + 5}^*(\alpha)},
$$
where $c_4 = c_3 e^{a_{\nu+3} + a_{\nu+4}}$.
%Take into account that for all $m_1, \ldots , m_n \in {\mathbb Z_+}$ 
%\begin{equation}
%(m_1 + \cdots + m_n)! \le e^{(n-1)(m_1 + \cdots + m_n)} m_1! \cdots  m_n! 
%\end{equation}
So if $k \in {\mathbb Z}_+$ then from the last inequality we get that
for all $\alpha \in {\mathbb Z_+^n}, \xi \in {\mathbb R}^n$
$$
(1 + \Vert \xi \Vert)^k \vert (D^{\alpha}f)(\xi) \vert  \le 
%M_{\nu, n, m} 
c_5 \Vert g \Vert_{k, \psi^*_{\nu}} \alpha! e^{-\psi_{\nu + 5}^*(\alpha)},
$$
where $c_5 > 0$ is some positive constant depending on $\nu, n$ and $k$.
By Theorem 2 $f$ can be holomorphically continued (uniquely) to entire function $F_f$ belonging to $E(\varPhi)$.   
Obviously by construction
$g = {\mathcal F}(F_f)$. 
The proof of Theorem 2 (see inequalities (2) and (7)) indicates that  there is a constant 
$c_6 > 0$  (depending on $\nu, n$ and $k$) such that for $z \in {\mathbb C}^n$
$$
(1 + \Vert z \Vert)^k \vert F_f(z) \vert \le c_6  
\Vert g \Vert_{k, \psi^*_{\nu}} e^{\varphi_{\nu + 8} (Im \ z_1 \vert, \cdots, \vert Im \ z_n \vert)}.
$$
Hence, 
$
p_{\nu + 8, k}(F_f) \le c_6 \Vert g \Vert_{k, \psi^*_{\nu}}. 
$
%где $K_{b, m}$ -- некоторая положительная постоянная, зависящая от $b$ и $m$. 
From this estimate it follows that the inverse mapping ${\mathcal F}^{-1}$ is continuous. 

Thus, we have proved that Fourier transform establishes a topological isomorphism between the spaces $E(\varPhi)$ and $G(\Psi^*)$. $\square$

\section{A special case of $\varPhi$}

{\bf 5.1. Some desired properties of the Young-Fenchel transform}. 
%In the proof of Theorem 4 Lemma 5 and Proposition 4 will be used. Also note that Proposition 4 plays an important role %in the proof of Proposition 3.

\begin{lemma} 
Let $u \in {\mathcal B}({\mathbb R}^n)$.
Then for each $\delta > 0$
$$
\displaystyle \lim_{x \to \infty} \frac {u^*((1+\delta) x) - u^*(x)}{\Vert x \Vert}= + \infty.
$$
\end{lemma} 

{\bf Proof}. Obviously, $u^*$ takes finite values on ${\mathbb R}^n$. 
For each $x \in {\mathbb R}^n$ denote by $\xi(x)$ a point where the supremum  of the function 
$g_x(\xi):= \langle x, \xi \rangle - u(\xi)$ over ${\mathbb R}^n$ is attained. 
Note that from the equality 
$u^*(x) + u(\xi(x))= \langle x,  \xi(x) \rangle$
and the fact that 
$
\displaystyle \lim_{x \to \infty} \frac {u^*(x)}{\Vert x \Vert}= + \infty
$
(arguments of Remark 2 can be applied here)
%(which easily follows from the inequality 
%$g^*(x) \ge \langle x, \xi \rangle - g(\xi)$ with $x, \xi \in {\mathbb R}^n$ by using the  continuity of $g$)
we have that 
$
\displaystyle \lim_{x \to  \infty} \frac {\langle x,  \xi(x) \rangle}{\Vert x \Vert}= + \infty
$
Now if $\delta > 0$ is arbitrary then from this and the inequality
$$
u^*((1+\delta) x)- u^*(x) \ge 
%(1+\delta) \langle x,  \xi(x) \rangle - g(\xi(x)) - \langle x,  \xi(x) \rangle + g(\xi(x))=
\delta \langle x,  \xi(x) \rangle, \ x \in {\mathbb R}^n,
$$
the assertion of lemma follows. $\square$

%\begin{proposition}
%Let $u \in {\mathcal A}_n \cap C^2({\mathbb R}^n)$ and be convex. Then 
%$$
%(u[e])^*(x) + (u^*[e])^*(x) = \sum \limits _{1 \le j \le n: x_j \ne 0} (x_j \ln x_j - x_j), \
%x \in [0, \infty)^n \setminus \{0\};
%$$
%$$
%(u[e])^*(0) + (u^*[e])^*(0) = 0.
%$$
%\end{proposition}

\begin{lemma}
Let $u \in {\mathcal B}({\mathbb R}^n)$.
Then
$$
(u[e])^*(x) + (u^*[e])^*(x) \le \sum \limits _{1 \le j \le n: \atop x_j \ne 0}
(x_j \ln x_j - x_j), \
x = (x_1, \ldots , x_n) \in [0, \infty)^n \setminus \{0\};
$$
$$
(u[e])^*(0) + (u^*[e])^*(0) \le 0.
$$
\end{lemma}

{\bf Proof}. For each $x =(x_1, \ldots , x_n) \in [0, \infty)^n$ and for each $\varepsilon > 0$ 
there are points  $t=(t_1, \ldots , t_n), \xi = (\xi_1, \ldots , \xi_n) \in {\mathbb R}^n$ such that 
$$
(u[e])^*(x) < \langle x, t \rangle - u[e](t) + \varepsilon, 
$$
$$
(u^*[e])^*(x) < \langle x, \xi \rangle - u^*[e](\xi) + \varepsilon.
$$
From this it follows that 
%Hence, 
%$$
%(u[e])^*(x) + (u^*[e])^*(x) < \langle x, t + \xi \rangle - u[e](t)  
%-\sup_{\eta \in {\mathbb R}^n} (\langle (e^{\xi_1}, \ldots , e^{\xi_n}),  \eta \rangle - u(\eta)) + 2 \varepsilon. 
%$$
%Thus, 
for each $\eta \in {\mathbb R}^n$ 
$$
(u[e])^*(x) + (u^*[e])^*(x) < \langle x, t + \xi \rangle - u[e](t) - \langle (e^{\xi_1}, \ldots , e^{\xi_n}),  \eta \rangle + u(\eta) + 2 \varepsilon. 
$$
Putting here $\eta = (e^{t_1}, \ldots , e^{t_n})$ we obtain that
$$
(u[e])^*(x) + (u^*[e])^*(x) < \sum \limits _{j=1}^n (x_j (t_j + \xi_j) - e^{\xi_j+ t_j})+ 2 \varepsilon. 
$$
Consequently,
$$
(u[e])^*(x) + (u^*[e])^*(x) < \sum \limits _{j=1}^n  
\sup_{y_j \in {\mathbb R}} (x_j y_j - e^{y_j})+ 2 \varepsilon.
$$
From this we get that
$$ 
(u[e])^*(x) + (u^*[e])^*(x) < \sum \limits _{1 \le j \le n: \atop x_j \ne 0} (x_j \ln x_j - x_j)+ 2 \varepsilon, \ x \in [0, \infty)^n \setminus \{0\};
$$
$$
(u[e])^*(0) + (u^*[e])^*(0) < 2 \varepsilon.
$$
Since $\varepsilon$ is arbitrary positive number 
then from the last two inequalities the assertion of Lemma follows. $\square$

\begin{proposition}
Let $u \in {\mathcal B}({\mathbb R}^n) \cap C^2({\mathbb R}^n)$ and be convex. Then 
$$
(u[e])^*(x) + (u^*[e])^*(x) = \sum \limits_{j =1}^n
%: x_j \ne 0} 
(x_j \ln x_j - x_j), \
x = (x_1, \ldots , x_n)\in (0, \infty)^n.
$$
\end{proposition}

%Proposition 3 follows from next two lemmas. 

%\begin{lemma}
%Let $u \in {\mathcal B}({\mathbb R}^n) \cap C^2({\mathbb R}^n)$ and be convex. Then 
%$$
%(u[e])^*(x) + (u^*[e])^*(x) \ge \sum \limits _{1 \le j \le n: x_j \ne 0} (x_j \ln x_j - x_j), \
%x \in [0, \infty)^n \setminus \{0\};
%$$
%$$
%(u[e])^*(0) + (u^*[e])^*(0) = 0.
%$$
%\end{lemma}

{\bf Proof}. Let $x = (x_1, \ldots , x_n) \in (0, \infty)^n$ be arbitrary. If we show that
$$
(u[e])^*(x) + (u^*[e])^*(x) \ge \sum \limits_{j =1}^n
%: \atop x_j \ne 0} 
(x_j \ln x_j - x_j),
$$
then (taking into account Lemma 6) the assertion will be proved.
%First assume that 
%Let $x = (x_1, \ldots , x_n) \in (0, \infty)^n$. 
First remark that for all $\xi, \mu \in {\mathbb R}^n$ 
$$
(u[e])^*(x) + (u^*[e] )^*(x) \ge \langle x,  \xi + \mu \rangle - (u[e]({\xi}) + u^*[e](\mu)).
$$
For an arbitrary $t = (t_1, \ldots , t_n) \in (0, \infty)^n $ denote $(\ln t_1, \ldots , \ln t_n)$ by $\xi(t)$ and
$(\ln \frac {x_1}{t_1}, \ldots , \ln \frac {x_n}{t_n})$ by $\mu(t)$ and put in the above inequality  
$\xi = \xi(t)$, $\mu = \mu(t)$. 
Then we get that
\begin{equation}
(u[e])^*(x) + (u^*[e])^*(x) \ge \sum \limits_{j=1}^n x_j \ln x_j - (u[e]({\xi(t)}) + u^*[e](\mu(t))).
\end{equation}

Further, note that there is a point $\zeta^* = (\zeta^*_1, \ldots , \zeta^*_n) \in {\mathbb R}^n$ where the supremum of the function 
$g_x: \zeta \in {\mathbb R}^n \to  \langle x, \zeta \rangle - u[e](\zeta)$ over ${\mathbb R}^n$
is attained.
Indeed, for each $\varepsilon > 0$ there exists a point  $\zeta(\varepsilon)=(\zeta_1(\varepsilon), \ldots , \zeta_n(\varepsilon)) \in {\mathbb R}^n$  such that 
\begin{equation}
(u[e])^*(x) < \langle x, \zeta(\varepsilon) \rangle - u[e](\zeta(\varepsilon)) + \varepsilon. 
\end{equation}
Remark that there exists a positive constant $C$ depending on $x$ such that 
$\Vert \zeta (\varepsilon) \Vert \le C$ for all $\varepsilon > 0$. Otherwise, 
there exists a decreasing to zero sequence 
$(\varepsilon_m)_{m=1}^{\infty}$ such that 
$\Vert \zeta (\varepsilon_m) \Vert \to + \infty$ as $m \to \infty$. 
From the growth conditions on $u$ we can find a constant $A>0$ such that
$$
u[e](\zeta) > e^{\zeta_1} + \cdots + e^{\zeta_n} - A, \ \zeta = (\zeta_1, \ldots , \zeta_n) \in {\mathbb R}^n.
$$ 
Then from this and the inequality (15) we obtain that  
$$
(u[e])^*(x) < \langle x, \zeta(\varepsilon_m) \rangle - e^{\zeta_1(\varepsilon_m)} - \cdots - e^{\zeta_n(\varepsilon_m)} + A + \varepsilon_m, \ m \in {\mathbb N}.
$$
From this inequality it follows that $(u[e])^*(x) = - \infty$. 
But it contradicts to the fact that 
$(u[e])^*(x) > - \infty$ (see Remark 2). 
Thus, we have shown that there exists a constant $C>0$ (depending on $x$) such that 
$\Vert \zeta_{\varepsilon} \Vert \le C$ for all $\varepsilon > 0$. 
Then using the Bolzano-Weierstrass theorem we can extract a sequence $(\zeta({\varepsilon_j}))_{j=1}^{\infty}$ converging to some point of ${\mathbb R}^n$. 
%(here $\lim \limits_{j \to \infty} {\varepsilon_j} = 0$). 
Denote this point by $\zeta^*$. Now from (15) we get that 
$
(u[e])^*(x) \le \langle x, \zeta^* \rangle - u[e](\zeta^*). 
$
From the other hand $(u[e])^*(x) \ge \langle x, \zeta \rangle - u[e](\zeta)$ for each $\zeta \in {\mathbb R}^n$. Hence, 
$(u[e])^*(x) = \langle x, \zeta^* \rangle - u[e](\zeta^*)$. 
Thus, $\zeta^*$ is the point where the supremum of the function $g_x$ over ${\mathbb R}^n$ is attained. 
Obviously, 
%since $u \in C^1({\mathbb R}^n)$ then 
$x_j = e^{\zeta^*_j} (D_j u)(e^{\zeta^*_1}, \ldots , e^{\zeta^*_n})$ ($j = 1, \ldots , n$).
Next, define a point $t^* = (t_1^*, \ldots , t_n^*) \in (0, \infty)^n $ by the rule $t_j^* = e^{\zeta^*_j}$. Then 
\begin{equation}
t_j^* (D_j u)(t^*) = x_j, \ j = 1, \ldots , n.
\end{equation}
%Let us show that $u[e]({\xi(t)}) + u^*[e](\mu(t)) \le \sum \limits _{j=1}^n x_j$.
Define the function 
$U_x: \eta \in {\mathbb R}^n \to u[e]({\xi(t^*)}) + \langle e^{\mu(t^*)}, \eta \rangle- u(\eta) - 
\sum \limits_{j=1}^n x_j$. 
In other words, 
$U_x(\eta)=u(t^*) + \sum \limits_{j=1}^n \frac {x_j \eta_j}{t_j^*} - u(\eta) - 
\sum \limits_{j=1}^n x_j, \ 
\eta \in {\mathbb R}^n$.
%If $\eta \ne t$ then by Taylor's formula we have that for some $\zeta \in {\mathbb R}^n$ (depending on $\eta$)   
%$$
%u(\eta) - u(t)=  \sum \limits_{\vert \alpha \vert =1} (D^{\alpha} u)(t)(\eta - t)^{\alpha} + 
%\sum \limits_{\vert \alpha \vert =2} (D^{\alpha} u)(\zeta)(\eta - t)^{\alpha}.
%$$ 
%$$u(t) - u(\eta) + \sum \limits_{j=1}^n \frac {x_j \eta_j}{t_j} - \sum \limits_{j=1}^n x_j = $$
If $\eta \ne t^*$ then using Taylor's formula we have that for some $\tau \in {\mathbb R}^n$ (depending on $\eta$) 
$$
U_x(\eta) = - \sum \limits_{\vert \alpha \vert =1} (D^{\alpha} u)(t^*)(\eta - t^*)^{\alpha} + \sum \limits_{j=1}^n x_j \frac {\eta_j - t_j^*}{t_j^*} -  \frac 1 2 \sum \limits_{\vert \alpha \vert =2} 
(D^{\alpha} u)(\tau)(\eta - t^*)^{\alpha}.
$$
Taking into account (16) we get that
$$
U_x(\eta) = - \frac 1 2 \sum \limits_{\vert \alpha \vert = 2} (D^{\alpha} u)(\tau)(\eta - t^*)^{\alpha}.
$$
Since $u$ is convex then from this equality 
it follows that $U_x(\eta) \le 0$. Also notice that $U_x(t^*) = 0$. 
Thus, $U_x(\eta) \le 0$ for all $\eta \in {\mathbb R}^n$. From this it follows that
$u[e]({\xi(t^*)}) + u^*[e](\mu(t^*)) \le \sum \limits _{j=1}^n x_j$. 
From the other hand for each $t \in (0, \infty)^n $ we have that
$$
u[e]({\xi(t)}) + u^*[e](\mu(t)) = u(t) + u^*\left(\frac {x_1}{t_1}, \ldots , \frac {x_n}{t_n}\right) 
\ge \sum \limits_{j=1}^n x_j.
$$  
Thus, $u[e]({\xi(t^*)}) + u^*[e](\mu(t^*)) = \sum \limits _{j=1}^n x_j$.
From this and (14) the desired inequality then follows. $\square$
%it follows that 
%$$
%(u[e])^*(x) + (u^*[e])^*(x) \ge \sum \limits_{j=1}^n (x_j \ln x_j - x_j), \
%x = (x_1, \ldots , x_n) \in (0, \infty)^n.
%$$

\begin{proposition}
Let $u \in {\mathcal A}({\mathbb R}^n) \cap C^2({\mathbb R}^n)$ and be convex. Then 
$$
(u[e])^*(x) + (u^*[e])^*(x) =  \sum \limits _{1 \le j \le n: \atop x_j \ne 0} (x_j \ln x_j - x_j), \,
x = (x_1, \ldots , x_n) \in [0, \infty)^n \setminus \{0\};
$$
$$
(u[e])^*(0) + (u^*[e])^*(0) = 0.
$$
\end{proposition}

{\bf Proof}. If $x \in (0, \infty)^n$ then the assertion follows from Proposition 3. 
Now let $x =(x_1, \ldots , x_n)$ belongs to the boundary of $[0, \infty)^n$ and $x \ne 0$. 
Assume for simplicity that the first $k$ ($1 \le k \le n-1$) coordinates of $x$ are positive and other coordinates are equal to zero. 
%For an arbitrary point $y = (y_1, \ldots , y_n) \in {\mathbb R}^n$ let $\tilde y = (y_1, \ldots , y_k)$. 
%Using these notations we have that
%According to our notations $x = (\tilde x, {\bf 0})$. 
%Denote points $(x_1, \ldots , x_k)$ by $\tilde x$, $(0, \ldots, 0) \in {\mathbb R}^{n-k}$ by ${\bf 0}$.  
%Then $x = (\tilde x, {\bf 0})$. 
For all $\xi = (\xi_1, \ldots , \xi_n), \mu = (\mu_1, \ldots , \mu_n) \in {\mathbb R}^n$ we have that
%\begin{equation}
%(u[e])^*(x) + (u^*[e])^*(x) \ge \langle \tilde x,  \tilde \xi + \tilde \mu \rangle - (u[e]({\xi}) + u^*[e](\mu)).
%\end{equation}
$$
(u[e])^*(x) + (u^*[e])^*(x) 
\ge \sum \limits _{j=1}^k x_j (\xi_j + \mu_j) - 
%\langle \tilde x,  \tilde \xi + \tilde \mu \rangle - 
(u(e^{\xi_1}, \ldots , e^{\xi_n}) + u^*(e^{\mu_1}, \ldots , e^{\mu_n})).
$$
From this inequality we get that  
$$
(u[e])^*(x) + (u^*[e])^*(x) 
\ge \sum \limits _{j=1}^k x_j (\xi_j + \mu_j) - 
$$
$$
-
(u(e^{\xi_1}, \ldots , e^{\xi_k}, 0, \ldots , 0) + u^*(e^{\mu_1}, \ldots , e^{\mu_k}, 0, \ldots , 0)).
$$
Let $\theta = (\theta_1, \ldots , \theta_k) \in (0, \infty)^k$ be arbitrary. 
Putting in the above inequality 
$\xi_j = \ln \theta_j$, $\mu_j = \ln \frac {x_j}{\theta_j}$ ($j = 1, \ldots , k)$ 
we obtain that
%Denote $(\ln \theta_1, \ldots , \ln \theta_k)$ by $\xi_k(t)$ and
%$(\ln \frac {x_1}{\theta_1}, \ldots , \ln \frac {x_k}{\theta_k})$ by $\mu_k(t)$. 
%Put in (20) $\xi = \xi_k(t)$ and $\mu = \mu_k(t)$. 
%Then we have 
\begin{equation}
(u[e])^*(x) + (u^*[e])^*(x) \ge \sum \limits_{j=1}^k x_j \ln x_j  
$$
$$
-
(u(\theta_1, \ldots , \theta_k, 0, \ldots , 0) + 
u^*(\frac {x_1}{\theta_1}, \ldots , \frac {x_k}{\theta_k}, 0, \ldots , 0)).
\end{equation}
Denote the point $(x_1, \ldots , x_k) \in {\mathbb R}^k $ by $\breve x$ and
define functions $u_k$ and ${\cal U}_{\breve x}$ on ${\mathbb R}^k$ by the rules: 
$$
u_k: \lambda = (\lambda_1, \ldots , \lambda_k) \in {\mathbb R}^k \to 
u(\lambda_1, \ldots , \lambda_k, 0, \ldots , 0);
$$
$$
{\cal U}_{\breve x}: \lambda = (\lambda_1, \ldots , \lambda_k) \in {\mathbb R}^k \to 
\langle \breve x,  \lambda \rangle -
u_k\left(e^{\lambda_1}, \ldots , e^{\lambda_k}\right).
$$
Repeating the same steps shown before in Proposition 3, we can find 
a point $\lambda^* = (\lambda_1^*, \ldots , \lambda_k^*) \in {\mathbb R}^k$ where the supremum 
of the function ${\cal U}_{\breve x}$ over ${\mathbb R}^k$ is attained. 
It is clear that 
$x_j = 
e^{\lambda^*_j} (D_j u_k)(e^{\lambda^*_1}, \ldots , 
e^{\lambda^*_k})$ ($j = 1, \ldots , k).$ 
Define a point $\theta^* = (\theta_1^*, \ldots , \theta_k^*) \in (0, \infty)^k$ by the rule $\theta_j^* = e^{\lambda^*_j}$ ($ j= 1, \ldots , k$). Then 
$
\theta_j^* (D_j u_k)(\theta^*) = x_j, \ j = 1, \ldots , k.
$
%Since $u_k$ is convex on ${\mathbb R}^k$, then 
Using similar computations as in Proposition 3 we obtain that 
$u_k(\theta_1^*, \ldots , \theta_k^*) + 
u_k^*\left(\frac {x_1}{\theta_1^*}, \ldots , \frac {x_k}{\theta_k^*}\right) = 
\sum \limits_{j=1}^k x_j$.
Now using that $u \in {\mathcal A}({\mathbb R}^n)$ we notice that
$$
u^*\left(\frac {x_1}{\theta_1}, \ldots , \frac {x_k}{\theta_k}, 0, \ldots , 0\right) = 
\sup \limits_{v = (v_1, \ldots , v_k, \ldots , v_n) \in {\mathbb R}^n} 
\left(\frac {x_1 v_1}{\theta_1} + \cdots + \frac {x_k v_k}{\theta_k} - u(v)\right) 
$$
$$
=
\sup \limits_{(v_1, \ldots , v_k) \in {\mathbb R}^k} 
\left(\frac {x_1 v_1}{\theta_1}  + \cdots + \frac {x_k v_k}{\theta_k} - u(v_1, \ldots , v_k, 0, \ldots , 0)\right) 
$$
$$
= u_k^*\left(\frac {x_1}{\theta_1}, \ldots , \frac {x_k}{\theta_k}\right).
$$
Thus, $u(\theta_1^*, \ldots , \theta_k^*, 0, \ldots , 0) + 
u^*\left(\frac {x_1}{\theta_1^*}, \ldots , \frac {x_k}{\theta_k^*}, 0, \ldots , 0 \right) = 
\sum \limits_{j=1}^k x_j$.
Finally, taking into account the inequality (17) we obtain that 
$$
(u[e])^*(x) + (u^*[e])^*(x) \ge \sum \limits_{j=1}^k (x_j \ln x_j - x_j).
$$
From this and the assertion of Lemma 6 the desired equality follows.

If $x =0$ then
$(u[e])^*(0) = - u(0)$, $(u^*[e])^*(0) = - \inf \limits_{\xi \in {\mathbb R}^n} u^*[e] (\xi) = -u^*(0) = u(0)$. Hence, 
$(u[e])^*(0) + (u^*[e])^*(0) = 0$. $\square$

%Thus, the proof of Lemma is complete.

\begin{corollary}
Let $u \in {\mathcal A}({\mathbb R}^n) \cap C^2({\mathbb R}^n)$ and be convex. Then 
$$
(u[e])^*(x) + (u^*[e])^*(x) \ge \sum \limits _{j=1}^n (x_j \ln (x_j  + 1) - x_j) - n, 
x = (x_1, \ldots , x_n) \in [0, \infty)^n.
$$
\end{corollary}

Notice that Propositions 3 and 4 are related to the following result obtained by S.V. Popenov 
(see Lemma 4 in \cite {N-P}): let $u \in {\mathcal A}({\mathbb R})$ be a convex function such that
$
\lim \limits_{x \to 0} \displaystyle \frac {u(x)} {x} = 0,
$
then
$$
(u[e])^*(x) + (u^*[e])^*(x) = x \ln x - x, \
x > 0,
$$
$$
(u(e))^*(0) + (u^*(e))^*(0) = 0.
$$

%In Proposition 3 we give an extension of this result for a case of 
%functions belonging to class ${\mathcal A}_n \cap C^2({\mathbb R}^n)$. But first let us prove the following two lemmas.

%From Lemmas 5 and 6 we have the following 

{\bf 5.2. Description of $E(\varPhi)$ by a system of weighted $C^{\infty}$-functions}. 
Choose a non-negative even function $\chi \in C_0^{\infty}({\mathbb R})$ with $supp \, \chi$ in $(-1, 1)$ and $\int_{\mathbb R} \chi(\xi) \ d \xi = 1$. 
Define a function $\omega$ on ${\mathbb R}^n$ by the rule: 
$\omega (x_1, \ldots , x_n)=\chi(x_1) \cdots \chi(x_n)$.
For each $m \in {\mathbb N}$ let 
$$
\varphi_{m, 1}(x) = \int_{{\mathbb R}^n} \varphi_m (x + \xi) \, \omega (\xi) \ d \lambda_n(\xi), \ x \in {\mathbb R}^n. 
$$
Here $d \lambda_n$ is the $n$-dimensional Lebesgue measure. 
The regularity properties of convolution ensures that $\varphi_{m, 1} \in C^{\infty}({\mathbb R}^n)$
and
$\varphi_{m, 1}(x_1, \dots , x_n) = \varphi_{m, 1}(\vert x_1 \vert, \dots , \vert x_n \vert)$ for $(x_1, \dots , x_n) \in {\mathbb R}^n$. Using convexity of $\varphi_m$ we have that
\begin{equation}
\varphi_m(x) \le \varphi_{m, 1}(x), \ x \in {\mathbb R}^n. 
\end{equation}
From this it follows that 
$\displaystyle \lim_{x \to \infty} \frac {\varphi_{m, 1}(x)}{\Vert x \Vert}= + \infty$. 
Using convexity of $\varphi_m$ and 
since  ${\varphi_m}_{|[0, \infty)^n}$ is not decreasing in each variable it is not difficult to show that ${\varphi_m}_{|[0, \infty)^n}$ is nondecreasing in each variable.
%Note that a restriction of $\varphi_{m, 1}$ to $[0, \infty)^n$ is nondecreasing in each variable. For example, show %that $\varphi_{m, 1}$ is nondecreasing in $x_1$ on $[0, \infty)$ while $x_2, \ldots , x_n$ are fixed and %non-negative. 
Thus, the family $\varPhi_1 = \{\varphi_{m, 1} \}_{m=1}^{\infty}$ is in ${\mathcal A}({\mathbb R}^n)$.

It is trivial to verify that for each $m \in {\mathbb N}$ and each $A > 0$ there exists a constant $s_{m, A} > 0$ such that 
\begin{equation}
\varphi_{m, 1}(x) + A \ln (1 + \Vert x \Vert) \le \varphi_{m+1, 1}(x) + s_{m, A}, \ x \in {\mathbb R}^n.
\end{equation}
%Thus, the family $\varPhi_1$ satisfies the condition of the form $i_0)$.

Further, 
%using the condition $i_2)$ on $\varPhi$  
for $x = (x_1, \ldots , x_n) \in [0, \infty)^n, \zeta = \zeta_1, \ldots , \zeta_n) \in [0, 1]^n$ we have that
$$
\varphi_{m, 1}(x + \zeta) = \int_{{\mathbb R}^n} \varphi_m (x + \zeta + \xi) \, \omega (\xi) \ d \lambda_n(\xi) 
$$
$$
=
\int_{{\mathbb R}^n} \varphi_m (\vert x_1 + \zeta_1 + \xi_1\vert, \ldots , \vert x_n + \zeta_n + \xi_n \vert) \, \omega (\xi) \ d \lambda_n(\xi).
$$
As ${\varphi_m}_{|[0, \infty)^n}$ is nondecreasing in each variable then 
$$
\varphi_{m, 1}(x + \zeta) \le \int_{{\mathbb R}^n} \varphi_m ((x_1 + \xi_1) + \vert \zeta_1\vert, \ldots , (x_n + \xi_n)  + \vert \zeta_n \vert) \, \omega (\xi) \ d \lambda_n(\xi).
$$
Now using the condition $i_2)$ on $\varPhi$ we have that
\begin{equation}
\varphi_{m, 1}(x + \zeta) \le \int_{{\mathbb R}^n} (\varphi_{m+1} (x_1 + \xi_1, \ldots , x_n + \xi_n) + K_m) \,
\omega (\xi) \ d \lambda_n(\xi)  
$$
$$
=\varphi_{m +1, 1}(x) + K_m.
\end{equation}
%Hence, the family $\varPhi_1$ satisfies the condition of the form $i_2)$.

Next, for $x = (x_1, \ldots , x_n) \in {\mathbb R}^n$ we have that
$$
\varphi_{m, 1}(2x) = 
%\int_{\mathbb R}^n \varphi_m (2 x + \xi) \, \omega (\xi) \ d \lambda_n(\xi)  
%$$
%$$
%=
%\int_{{\mathbb R}^n} \varphi_m (2x_1 + \xi_1, \ldots , 2x_n + \xi_n) \, \omega (\xi) \ d \lambda_n(\xi)  
%$$
%$$
%=
\int_{{\mathbb R}^n} \varphi_m (\vert 2x_1 + \xi_1 \vert, \ldots , \vert 2x_n + \xi_n \vert) \, \omega (\xi) \ d \lambda_n(\xi).
$$
From this using nondecreasity in each variable of ${\varphi_m}_{|[0, \infty)^n}$ we have that
$$
\varphi_{m, 1}(2x) \le \int_{{\mathbb R}^n} \varphi_m (\vert 2x_1 \vert + \vert \xi_1 \vert, \ldots , \vert 2x_n \vert + \vert \xi_n \vert)\, \omega (\xi) \ d \lambda_n(\xi).
$$
Now due to the condition $i_2)$ on $\varPhi$ we have that
$$
\varphi_{m, 1}(2x) \le \int_{{\mathbb R}^n} 
(\varphi_{m+1} (\vert 2x_1 \vert, \ldots , \vert 2x_n \vert) + K_m) \, \omega (\xi) \ d \lambda_n(\xi) = 
\varphi_{m+1} (2x) + K_m.
$$
Thanks to the condition $i_3)$ on $\varPhi$ we get that
\begin{equation}
\varphi_{m, 1}(2x) \le  
\varphi_{m+2} (x) + K_m + a_{m+1}.
\end{equation}
%Earlier it was noticed that for each $k \in {\mathbb N}$ 
%$\varphi_k(x) \le \varphi_{k, 1}(x), \ x \in {\mathbb R}^n$.
Using the inequality (18) we obtain that
\begin{equation}
\varphi_{m, 1}(2x) \le  
\varphi_{m+2, 1} (x) + K_m + a_{m+1}, \ x \in {\mathbb R}^n.
\end{equation}

Now introduce the family 
$\varPhi_2 =\{\varphi_{2m, 1}\}_{m=1}^{\infty}$. Obviously, $\varPhi_2$ is in ${\mathcal A}({\mathbb R}^n)$. From the inequalities (19), (20) and (22) it follows that $\varPhi_2$ satisfies the condition of the form $i_0)$, $i_2)$ and $i_3)$. 

From the inequality (21) we have that 
$$
\varphi_{2m, 1}(x) \le  \varphi_{2m+2} (x) + K_{2m} + a_{2m+1}, \ x \in {\mathbb R}^n.
$$
This inequality and the inequality (18) mean that $E(\varPhi)$ can be described by the system $\varPhi_2$.

{\bf 5.3. Proof of Theorem 4}. 
We may assume that functions of the family $\varPhi$ belong to $C^{\infty}({\mathbb R}^n)$ (see subsection 5.2). 
Further, note that using convexity of functions of the family $\varPhi$ and the condition $i_3)$ 
on $\varPhi$ we have that for each $k \in {\mathbb N}$ 
$$
2 \varphi_k(x) \le \varphi_{k+1}(x) + a_k + \varphi_k(0), \ x \in {\mathbb R}^n. 
$$
This means that the condition $i_4)$ holds in our case with 
$l_k = a_k + \varphi_k(0)$. 
From the last inequality it follows that for each $k \in {\mathbb N}$ we have that 
$$
2 \psi_k(x) \le \psi_{k+1}(x) + a_k + \varphi_k(0), \ x \in {\mathbb R}^n. 
$$
Hence, the inequality (12) holds in our case. 
%Recall that by this ineaquality for each $k \in {\mathbb N}$ we have that
%$$
%\psi_{k+1}^*(x+y) \le \psi_k^*(x) + \psi_k^*(y) + l_k, \ x, y \in [0, \infty)^n.  
%$$

Now let $\nu \in {\mathbb N}$ and $f \in G(\psi_{\nu}^*)$. 
Let $m \in {\mathbb Z}_+$ be arbitrary. 
%Let $x \in {\mathbb R}^n$ be arbitrary. 
For all $\alpha \in {\mathbb Z_+^n}$ with $\vert \alpha \vert \le m$,  
$\beta  \in {\mathbb Z}_+^n$ and 
$x = (x_1, \ldots , x_n) \in {\mathbb R}^n$ with non-zero coordinates we have that
$$
\vert (D^{\alpha} f)(x) \vert \le 
\frac 
{\Vert f \Vert_{m, \psi_{\nu}^*}
\beta! 
e^{-\psi_{\nu}^*(\beta)}}
{\prod \limits_{j=1}^n \vert x_j \vert^{\beta_j}} \ .
$$
Take into account that 
$j! \le \frac {(j+1)^{j+1}}{e^j}$ for all $j \in {\mathbb Z_+}$. 
Then
\begin{equation}
\vert (D^{\alpha} f)(x) \vert \le \Vert f \Vert_{m, \psi_{\nu}^*}  
e^{-\psi_{\nu}^*(\beta)} 
{\prod \limits_{j=1}^n \frac {(\beta_j + 1)^{\beta_j + 1}}
{(e\vert x_j \vert) \vert^{\beta_j}}} \ .
\end{equation}
Our aim is to obtain a suitable estimate of 
$e^{-\psi_{\nu}^*(\beta)} 
{\prod \limits_{j=1}^n \frac {(\beta_j + 1)^{\beta_j + 1}}
{(e\vert x_j \vert) \vert^{\beta_j}}}$ from above.
For $\beta \in {\mathbb Z}_+^n$ let 
$\Omega_{\beta} = \{\xi = (\xi_1, \ldots , \xi_n) \in {\mathbb R}^n: 
\beta_j \le \xi_j  < \beta_j + 1 \ (j =1, \ldots , n) \}$. 
Also, for $\lambda > 0$ let $\tilde \lambda:  = \max (\lambda, 1)$.
Using the inequality (12) and nondecreasity in each variable of $\psi^*$  in $[0, \infty)^n$ 
we have that for $\xi \in \Omega_{\beta}$ and $\mu = (\mu_1, \ldots , \mu_n) \in (0, \infty)^n $  
$$
e^{-\psi_{\nu}^*(\beta)} 
{\prod \limits_{j=1}^n \frac {(\beta_j + 1)^{\beta_j + 1}}
{\mu_j^{\beta_j}}} 
\le e^{-\psi_{\nu + 1}^*(\xi) +\psi_{\nu}^*(1, \ldots , 1) + l_{\nu}}
{\prod \limits_{j=1}^n \frac {\tilde \mu_j (\xi_j + 1)^{\xi_j + 1}}
{\mu_j^{t_j}}}.
$$
%Recall that by our notations $\tilde \mu_j = \max (1, \mu_j)$. 
Denoting $e^{\psi_{\nu}^*(1, \ldots , 1) + l_{\nu}}$ by $C_1$ we 
rewrite the last inequality in the following form
$$
e^{-\psi_{\nu}^*(\beta)} 
{\prod \limits_{j=1}^n \frac {(\beta_j + 1)^{\beta_j + 1}}
{\mu_j^{\beta_j}}} 
\le C_1 
%e^{-\psi_{\nu + 3n + 1}^*(t)} 
e^{\sum \limits_{j=1}^n (\ln \tilde \mu_j  + 
(t_j + 1) \ln (t_j + 1) - t_j \ln \mu_j) -\psi_{\nu + 1}^*(t)}.
$$
Now using the Corollary 2  we obtain that
$$
e^{-\psi_{\nu}^*(\beta)} 
{\prod \limits_{j=1}^n \frac {(\beta_j + 1)^{\beta_j + 1}}
{\mu_j^{\beta_j}}} 
\le C_2 
e^{\sum \limits_{j=1}^n (\ln \tilde \mu_j  + \ln (t_j + 1) - t_j \ln \mu_j + t_j) + 
(\varphi_{\nu + 1}^*[e])^*(t)},
$$
where $C_2 = C_1 e^n$.
Obviously, there exists a constant $C_3 > 0$ such that 
$$
e^{-\psi_{\nu}^*(\beta)} 
{\prod \limits_{j=1}^n \frac {(\beta_j + 1)^{\beta_j + 1}}
{\mu_j^{\beta_j}}}  \le 
C_3 
e^{(\varphi_{\nu + 1}^*[e])^*(t) - \sum \limits_{j=1}^n  t_j \ln \frac {\mu_j}{4} + 
\sum \limits_{j=1}^n \ln \tilde \mu_j}.
$$
From this it follows that
$$
\inf_{\beta \in {\mathbb Z}_+^n}
e^{-\psi_{\nu}^*(\beta)} 
{\prod \limits_{j=1}^n \frac {(\beta_j + 1)^{\beta_j + 1}}
{\mu_j^{\beta_j}}}  \le 
C_3 
e^{\inf \limits_{t \in [0, \infty)^n} ((\varphi_{\nu + 1}^*[e])^*(t) - \sum \limits_{j=1}^n  t_j \ln \frac {\mu_j}{4}) + 
\sum \limits_{j=1}^n \ln \tilde \mu_j}.
$$
In other words,
$$
\inf_{\beta \in {\mathbb Z}_+^n}
e^{-\psi_{\nu + 3n}^*(\beta)} 
{\prod \limits_{j=1}^n \frac {(\beta_j + 1)^{\beta_j + 1}}
{\mu_j^{\beta_j}}}  \le 
C_3 
e^{-\sup \limits_{t \in [0, \infty)^n} (\sum \limits_{j=1}^n  t_j \ln \frac {\mu_j}{4} - 
(\varphi_{\nu + 1}^*[e])^*(t)) + 
\sum \limits_{j=1}^n \ln \tilde \mu_j}.
$$
Taking into account, by Remark 2, that the function  
$(\varphi_{\nu + 1}^*[e])^*$ takes finite values on $[0, \infty)^n$ and $(\varphi_{\nu + 1}^*[e])^*(x) = +\infty$ if 
$x \notin [0, \infty)^n$ we can rewrite the last inequality in the following form
$$
\inf_{\beta \in {\mathbb Z}_+^n}
e^{-\psi_{\nu + 3n}^*(\beta)} 
{\prod \limits_{j=1}^n \frac {(\beta_j + 1)^{\beta_j + 1}}
{\mu_j^{\beta_j}}}  \le 
C_3 
e^{-\sup \limits_{t \in {\mathbb R}^n} (\sum \limits_{j=1}^n  t_j \ln \frac {\mu_j}{4} - 
(\varphi_{\nu + 1}^*[e])^*(t)) + 
\sum \limits_{j=1}^n \ln \tilde \mu_j}.
$$
Note that by Lemma 4 the function $\varphi_{\nu + 1}^*[e]$ is convex on ${\mathbb R}^n$ with finite values (thus, $\varphi_{\nu + 1}^*[e]$ is continuous on ${\mathbb R}^n$ (see \cite {R}, Corollary 10.1.1)).
Taking into account that the Young-Fenchel comjugation is involutive (see \cite {R}, Theorem 12.2) we have that 
$$
\sup \limits_{t \in [0, \infty)^n} (\sum \limits_{j=1}^n  t_j \ln \frac {\mu_j}{4} - 
(\varphi_{\nu + 1}^*[e])^*(t)) = 
\varphi_{\nu + 1}^*[e]
\left(\ln \frac {\mu_1}{4}, \ldots , \ln \frac {\mu_n}{4}\right) .
$$
Thus,
$$
\inf_{\beta \in {\mathbb Z}_+^n}
e^{-\psi_{\nu}^*(\beta)} 
{\prod \limits_{j=1}^n \frac {(\beta_j + 1)^{\beta_j + 1}}
{\mu_j^{\beta_j}}}  \le 
C_3 
e^{-\varphi_{\nu + 1}^*[e]
\left(\ln \frac {\mu_1}{4}, \ldots , \ln \frac {\mu_n}{4}\right) + 
\sum \limits_{j=1}^n \ln \tilde \mu_j}.
$$
In other words, 
\begin{equation}
\inf_{\beta \in {\mathbb Z}_+^n}
e^{-\psi_{\nu}^*(\beta)} 
{\prod \limits_{j=1}^n \frac {(\beta_j + 1)^{\beta_j + 1}}
{\mu_j^{\beta_j}}}  \le 
C_3 
e^{-\varphi_{\nu + 1}^*(\frac {\mu}{4}) + 
\sum \limits_{j=1}^n \ln \tilde \mu_j}.
\end{equation}
Note that using the condition $i_3)$ on $\varPhi$ it is easy to obtain that for each $j \in {\mathbb N}$
\begin{equation}
\varphi_{j+1}^*(\xi) \le \varphi_j^*\left(\frac {\xi}{2}\right) + a_j, \ \xi \in {\mathbb R}^n.
\end{equation}
Due to this inequality we get from (24) that 
\begin{equation}
\inf_{\beta \in {\mathbb Z}_+^n}
e^{-\psi_{\nu}^*(\beta)} 
{\prod \limits_{j=1}^n \frac {(\beta_j + 1)^{\beta_j + 1}}
{\mu_j^{\beta_j}}}  \le 
C_4 
e^{-\varphi_{\nu + 3}^*(\mu) + 
\sum \limits_{j=1}^n \ln \tilde \mu_j},
\end{equation}
where $C_4 = C_3 e^{a_{\nu + 1} + a_{\nu + 2}}$.
Also note that from (25) it follows that for each $j \in {\mathbb N}$
$$
\varphi_j^*(\xi) - \varphi_{j+1}^*(\xi) \ge \varphi_j^*(\xi) - 
\varphi_j^*\left(\frac {\xi}{2}\right) - a_j, \ \xi \in {\mathbb R}^n.
$$
From this and Lemma 5 we get that
$$
\displaystyle \lim_{\xi \to \infty} \frac 
{\varphi_j^*(\xi) - \varphi_{j+1}^*(\xi)}{\Vert \xi \Vert}= + \infty.
$$
Using this we obtain from (26) that
$$
\inf_{\beta \in {\mathbb Z}_+^n}
e^{-\psi_{\nu}^*(\beta)} 
{\prod \limits_{j=1}^n \frac {(\beta_j + 1)^{\beta_j + 1}}
{\mu_j^{\beta_j}}}  \le 
C_5 
e^{-\varphi_{\nu + 4}^*(\mu)},
$$
where
$C_5$ is some positive number depending on $\nu$.
From this and the inequality (23) we obtain at last that for all 
$x = (x_1, \ldots , x_n) \in {\mathbb R}^n$ with non-zero coordinates and for all 
$\alpha \in {\mathbb Z_+^n}$ with $\vert \alpha \vert \le m$
$$
\vert (D^{\alpha} f)(x) \vert \le  C_5 \Vert f \Vert_{m, \psi_{\nu}^*}
e^{-\varphi_{\nu + 4}^*(e \vert x_1 \vert, \ldots , e \vert x_n \vert)} \le
C_5 \Vert f \Vert_{m, \psi_{\nu}^*}
e^{-\varphi_{\nu + 4}^*(x)}.
$$
%where $C_6 = 3^n C_5$.

Obviously, the last inequality holds for all $x \in {\mathbb R}^n$.  
%Since $\varphi_{\nu + 3n + 2}^*$ is nondecreasing on $\overline \Pi$ in each variable then
%$$
%\vert (D^{\alpha} f)(x) \vert \le  C_6 N_{\nu + n, m}(f) 
%e^{-\varphi_{\nu + 3n + 2}^*(\vert x_1 \vert, \ldots , \vert x_n \vert)},
%$$
%where $C_6 = C_0 C_5$.
Thus,  
$$
q_{m, \nu + 4}(f) \le C_5 \Vert f \Vert_{m, \psi_{\nu}^*}, \ f \in G(\psi_{\nu}^*).
$$
From this it follows that the identity mapping $J$ acts from $G(\Psi^*)$ to $GS(\varPhi^*)$ and is continuous.

We proceed to show that $J$ is surjective. Let $f \in GS(\varPhi^*)$. 
Then $f \in GS(\varphi_{\nu}^*)$ for some $\nu \in {\mathbb N}$. 
Fix $m \in {\mathbb Z}_+$. 
Consider an arbitrary point $x = (x_1, \ldots , x_n) \in {\mathbb R}^n$  with non-zero coordinates.
For all $\alpha \in {\mathbb Z_+^n}$ with $\vert \alpha \vert \le m$ we have that
$$
\vert (D^{\alpha}f)(x) \vert \le q_{m, \nu}(f)
e^{-\varphi_{\nu}^*(\vert x_1 \vert, \ldots , \vert x_n \vert)}.
$$
In other words,
$$
\vert (D^{\alpha}f)(x) \vert   \le q_{m, \nu}(f)
e^{-\varphi_{\nu}^*[e](\ln \vert x_1 \vert, \ldots , \ln \vert x_n \vert)}.
$$
From this we have that
$$
\vert (D^{\alpha}f)(x) \vert  \le q_{m, \nu}(f)  
\exp (- \sup \limits_{t =(t_1, \ldots , t_n) \in [0, \infty)^n} 
(\sum \limits_{j=1}^n t_j \ln \vert x_j \vert - (\varphi_{\nu}^*[e])^*(t))).
$$
Now using Proposition 4 we get that
$$
\vert (D^{\alpha}f)(x) \vert \le q_{m, \nu}(f)  
e^{- \sup \limits_{t =(t_1, \ldots , t_n) \in [0, \infty)^n} 
(\sum \limits_{1 \le j \le n:  t_j \ne 0} (t_j \ln (e \vert x_j \vert)  - t_j \ln t_j)+ 
\psi_{\nu}^*(t))}.
$$
Consequently, if $\beta \in {\mathbb Z}_+^n$ then
$$
\vert (D^{\alpha}f)(x) \vert \le q_{m, \nu}(f) e^{-\psi_{\nu}^*(\beta)}
\prod \limits_{1 \le j \le n: \beta_j \ne 0} 
\frac 
{\beta_j^{\beta_j}}
{(e \vert x_j \vert)^{\beta_j}} 
$$
From this we finally obtain that for all $x \in {\mathbb R}^n$  with non-zero coordinates
$$
\vert x^{\beta}(D^{\alpha}f)(x) \vert  \le 
q_{m, \nu}(f) 
\beta!  e^{-\psi_{\nu}^*(\beta)}, \ \beta \in {\mathbb Z}_+^n, \vert \alpha \vert \le m.
$$
Clearly, this inequality holds for all $x \in {\mathbb R}^n$. 
Thus,
$$
\Vert f \Vert_{m, \psi_{\nu}^*} \le q_{m, \nu}(f).
$$
Since here $m \in {\mathbb Z}_+$ is arbitrary then $f \in G(\psi_{\nu}^*)$. Hence, $f \in G(\Psi^*)$.
Also from the last inequality it follows that the mapping $J^{-1}$ is continuous. Thus, 
the spaces $G({\Psi^*})$ and $GS({\varPhi^*})$ coincide. $\square$

{\bf Acknowledgements}. The author is very grateful to the referee for careful reading, valuable comments and suggestions. The research was supported by grants from RFBR (15-01-01661). 

\pagebreak


\begin{thebibliography}{99}

\bibitem {B-R}
 J.J. Betancor, L. Rodr{\'\i}guez-Mesa, Characterizations of $W$-type spaces, {\it Proc. Amer. Math. Soc.}, 
{\bf 126} (1998), no. 5, 1371--1379.

\bibitem {Br} 
N.G. De Bruijn, A theory of generalized functions, with applications to Wigner distribution and Weyl correspondence, {\it Nieuw Archief Wiskunde}, {\bf 21} (1973), 205--280.

\bibitem {Cho}
J. Cho, A characterization of Gelfand-Shilov space based on Winger distribution, {\it Comm. Korean Math. Soc.}, {\bf 14} no. 4 (1999), 761--767.

\bibitem {C-C-K 1}  J. Chung, S.-Y. Chung, and D. Kim,
Characterizations of the Gelfand-Shilov spaces via Fourier transforms, {\it Proc. Amer. Math. Soc.}, {\bf 124} (1996), 7, 2101--2108.

\bibitem {C-C-K 2} J. Chung, S.-Y. Chung, and D. Kim, Equivalence of the Gelfand-Shilov spaces, {\it Journal of Math. Anal. and Appl.}, {\bf 203} (1996), 828--839.

\bibitem {E-K} 
S.J.L. van Eijndhoven, M.J. Kerkhof, The Hankel transformation and spaces of type W, {\it Reports on Applied and Numerical Analysis}, Departament of Mathematics and Computing Science, Eindhoven University of Technology, 1988.

\bibitem {FS}
V.Ya. Fainberg, M.A. Soloviev, Nonlocalizability and asymptotical commutativity, {\it Theoretical and Mathematical Physics}, {\bf 93} no. 3 (1992), 1438--1449.

\bibitem {GS1} 
I.M. Gelfand, G.E. Shilov, {\it Generalized functions}, Vol. 2, Academic Press, New York,
1968.

\bibitem {GS2} 
I.M. Gelfand, G.E. Shilov, {\it Generalized functions}, Vol. 3, Academic Press, New York, 1967. 

\bibitem {Gur1} 
B.L. Gurevich, New types of spaces of fundamental and generalized
functions and Cauchy's problem for systems of finite difference equations,
{\it Doklady Akad. Nauk SSSR (N.S.)}, {\bf 99} (1954), 893--895.

\bibitem {Gur2}
B.L. Gurevich, New types of fundamental and generalized spaces and
Cauchy's problem for systems of difference equations involving differential
operations, {\it Doklady Akad. Nauk SSSR (N.S.)}, {\bf 108} (1956), 1001--1003.

\bibitem {J-E} 
A.J.E.M. Janssen, S.J.L. van Eijndhoven, Spaces of type W, growth of Hermite coefficients, Wigner distribution and Bargmann transform, 
{\it J. Math. Anal. Appl.}, {\bf 152} no. 2 (1990), 368--390.

%\bibitem {K-R} Krasnosel'skii M.A., Rutickii Ya. B., {\it Convex functions and Orlicz spaces}. P. Noordhoff LTD, Groningen, 1961.

%\bibitem {M-P} 
%Musin I.Kh., Popenov S.V., {\it On a weighted space of infinitely differentiable functions in ${\mathbb R}^n$}, Ufa Mathematical Journal, {\bf 2} %(2010), 54--62. 

%\bibitem {M} 
%Musin I.Kh., {\it Approximation by polynomials in a weighted space of infinitely differentiable functions with an application to hypercyclicity}. %Extracta Mathematicae, {\bf 27} (2012), 75--90. 

\bibitem{N-P}
V.V. Napalkov, S.V. Popenov, On Laplace transformation on weighted Bergman space of entire functions 
on ${\mathbb C}^n$. {\it Doklady Mathematics},  {\bf 55}, no. 1 (1997), 110--112.

\bibitem {P}
R.S. Pathak, On Hankel transformable spaces and a Cauchy problem, {\it Can. J. Math.}, 
{\bf 37} (1985), 84--106.

\bibitem {P-U} R.S. Pathak, S.K. Upadhyay, $W^p$-spaces and Fourier transform, 
{\it Proc. Amer. Math. Soc.}, {\bf 121} no. 3 (1994), 733--738.

\bibitem {R}
R.T. Rockafellar, {\it Convex analysis},  Princeton, New Jersey, Princeton University Press, 1970.

\bibitem {U} S.K. Upadhyay, W-Spaces and Pseudo-differential Operators, {\it Applicable Analysis}, {\bf 82} no. 4 (2003), 381--397.

\end{thebibliography}
\end{document}